\documentclass[12pt]{article}
\usepackage{amsfonts}
\usepackage{amssymb}
\usepackage{polski}

\input{tcilatex}
\begin{document}

\begin{center}
\textbf{Sequential regular variation: extensions of Kendall's theorem}\\[0pt]
\medskip
\noindent\textbf{ by }\\[0pt]

\bigskip

\textbf{N. H. Bingham and A. J. Ostaszewski}\\[0pt]

\bigskip

\textit{David George Kendall 1918-2007, in memoriam}

\textit{To John Kingman, on his 80}$^{\text{\textit{th}}}$ \textit{birthday
(28 August 1939)}

\bigskip
\end{center}

\noindent \textbf{Abstract. }Regular variation is a continuous-parameter
theory; we work in a general setting, containing the existing Karamata,
Bojanic-Karamata/de Haan and Beurling theories as special cases. We give
sequential versions of the main theorems, that is, with sequential rather
than continuous limits. This extends the main result, a theorem of Kendall's
(which builds on earlier work of Kingman and Croft), to the general setting.

\bigskip

\noindent \textbf{Keywords. }Kendall's Theorem, regular variation, quasi
regular variation, general regular variation, uniform convergence theorem, Go%
\l \k{a}b-Schinzel equation\textit{,} Beurling-Goldie equation, essential
limits, croftian theorems, category-measure duality.

\bigskip

\noindent \textbf{Classification}: 26A03, 26A12, 33B99, 39B22.

\bigskip

\bigskip

\textbf{1. Introduction. }

We are concerned here with \textit{regular variation}, RV\ for short (for
background on which we refer to the standard work [BinGT], BGT below for
brevity). This occurs first in the classical \textit{Karamata }setting%
\begin{equation}
\lim\nolimits_{x\rightarrow \infty }f(tx)/f(x)=K(t)\qquad (\text{locally
uniformly in }t\text{ (}\forall t>0)\text{)}  \tag{$K$}
\end{equation}%
(see e.g. BGT Ch. 1), and then in the \textit{Bojanic-Karamata/de Haan }%
setting (see e.g. BGT Ch. 3),%
\begin{equation}
\lbrack f(x+t)-f(x)]/h(x)\rightarrow K(t)\qquad (\text{locally uniformly in }%
t)\text{.}  \tag{$BKdH$}
\end{equation}%
The basic result in the theory of regular variation is the\textit{\ Uniform
Convergence Theorem}, UCT (BGT, Th. 1.2.1): with $f$ and $h$ above Baire
(i.e. having the Baire property, BP) or measurable, convergence is
necessarily locally uniform, hence our general assumption both here and in
the contexts below.

Next, we need the \textit{Beurling }setting (see e.g. BGT \S\ 2.11 and
[BinO5, 10]), using an auxiliary function $\varphi ,$ best implemented
algebraically (see \S\ 3 below) via the \textit{Popa binary operation}
\[
x\circ _{\varphi }t:=x+t\varphi (x).
\]%
This leads to the \textit{Beurling} theory of regular variation [BinO10],
below.

The auxiliary $\varphi $ above is to satisfy%
\[
\eta _{x}(t):=\varphi (x+t\varphi (x))/\varphi (x)=\varphi (x\circ _{\varphi
}t)/\varphi (x)\rightarrow \eta (t)\qquad (\text{locally uniformly in }t)%
\text{.}
\]%
When $\varphi (x)=O(x),$ these are the \textit{self-equivarying} functions, $%
\varphi \in SE,$ of [Ost2]; when specialized to the case $\varphi (x)=o(x)$
and $\eta \equiv 1$ these reduce to the classical \textit{self-neglecting
functions} [BGT, \S\ 2.11]. Here again, it is known that, for $\varphi $
Baire/measurable, the convergence is necessarily locally uniform, provided $%
\varphi $ satisfies one of a series of four possible additional properties,
including continuity or monotonicity: see [Ost2, Th. 4]. We will need

\bigskip

\noindent \textbf{Theorem O }[Ost2, Th. 0].\textbf{\ }\textit{If }$\varphi
(x)=O(x)$\textit{\ and }$\eta _{x}(t)\rightarrow \eta (t)=\eta ^{\varphi
}(t),$ \textit{locally uniformly in }$t$\textit{, then }$\eta $\textit{\
satisfies the Go\l \k{a}b-Schinzel functional equation}%
\begin{equation}
\eta (s\circ _{\eta }t)=\eta (s)\eta (t).  \tag{$GS$}
\end{equation}

\noindent \textit{Notational convention. }In Theorem O above $\eta _{x}$
contains the $x$ which tends to infinity. After this passage to the limit,
attention focuses on the limit function $\eta (t)$ which will depend on a
parameter $\rho $, below. We allow ourselves to denote this limit by $\eta
_{\rho }(t)$ and let context speak for itself here.

\bigskip

For $\varphi $ Baire/measurable $\eta ^{\varphi }$ is measurable and so
continuous, as are the positive solutions of ($GS$), the only ones of
interest here, which take the form $\eta (t)=\eta _{\rho }(t):=1+\rho t,$
for $t>\rho ^{\ast }:=-1/\rho $ with $\rho \geq 0$ (and $0$ to the left of $%
\rho ^{\ast },$ though here we work in $\mathbb{R}_{+}$) -- see the surveys
[Brz] and [Jab5]; cf. [Ost2]. So $\eta =\eta _{0}\equiv 1$ yields the
desired limit of $\eta _{x}(t)$ for the self-neglecting case.

Finally, we need the \textit{general }setting of our title,
\begin{equation}
\lbrack f(x\circ _{\varphi }t)-f(x)]/h(x)\rightarrow K(t)\qquad (\text{%
locally uniformly in }t),  \tag{$GRV$}
\end{equation}%
recently developed in [BinO10]. Here $f$ is the function of primary
interest; $h$ and $\varphi $ are auxiliary functions; the limit function on
the right we call the \textit{kernel} function. The instance here $%
h(x)\equiv 1$ pre-dating this, developed in [BinO5] is termed Beurling RV, $%
f $ then is said to be $\varphi $-RV: again see \S\ 3.

These (limit) kernel functions $K$ satisfy \textit{functional equations.}
The classical Karamata setting yields the multiplicative\textit{\ Cauchy
functional equation, }($CFE$) for short below.\textit{\ }In addition to the
\textit{Go\l \k{a}b-Schinzel }functional equation above,\textit{\ }there
are: the \textit{Chudziak-Jab\l o\'{n}ska }functional equation%
\begin{equation}
K(u\circ _{\eta }v)=K(v)K(u);  \tag{$CJ$}
\end{equation}%
the \textit{Beurling-Goldie }functional equation of the general setting%
\begin{equation}
K(u\circ _{\eta }v)=K(u)\circ _{\sigma }K(v),  \tag{$BG$}
\end{equation}%
and the original \textit{Goldie }functional equation [Ost3] of the $h\equiv
1 $ setting, which it extends,%
\begin{equation}
K(u+v)=K(u)\circ _{\sigma }K(v).  \tag{$G$}
\end{equation}%
See \S\ 9.6 for their continuous solutions, and a discussion of how
discontinuous solutions are excluded by the `blanket assumption of
non-triviality', stated in \S ~3 ahead of Theorem 4 (after the necessary
preliminaries in \S ~2).

All four of the limiting settings above involve \textit{continuous }limits.
However, \textit{sequential }limits (see e.g. BGT \S\ 1.9) are also
important, both in theory (see the theorems below) and in applications
(particularly to probability -- see e.g. \S\ 9.3 -- which as it happens
originally motivated the theory).

The prototypical sequential result here is due to Croft [Cro] in 1957. The
role of the Baire Category Theorem, and the relevance to probability theory,
are due to Kingman [Kin1,2] in 1963 and 1964, and Kendall [Ken] in 1968.

The Baire Category Theorem is sequential, and so its role in the sequential
results here is thematic. All the `Baire' results below need only the Axiom
of Dependent Choice(s), DC. We comment briefly on the set-theoretic
axiomatic background here in \S\ 9.1.

For the interplay between category and measure in settings such as this, we
refer to a number of our previous studies, for instance [BinO5,6,7] and
[Ost1,2]; it is \textit{category rather than measure} that is primary here.
For background on the axiomatics underpinning results in this area, we refer
to our recent survey [BinO9]; see again \S\ 9.1.

We will rely on the following combinatorial tool. Below, $B$ is
`negligible', $B\in \mathcal{N}$, will mean $B$ is meagre or null according
to context, `quasi all' will mean `off a negligible set', while
`non-negligible' will implicitly mean Baire/measurable (and
non-meagre/non-null).

\bigskip

\noindent \textbf{Proposition 1 (Affine Two-sets Lemma }[BinO5, Lemma 2]%
\textbf{). }\textit{For }$c_{n}\rightarrow c>0$ \textit{and }$%
z_{n}\rightarrow 0,$ \textit{if }$cB\subseteq A$ \textit{for }$A,B$ \textit{%
non-negligible, then for quasi all }$b\in B$ \textit{there exists an
infinite set }$\mathbb{M}=\mathbb{M}_{b}\subseteq \mathbb{N}$ \textit{such
that}%
\[
\{c_{m}b+z_{m}:m\in \mathbb{M}\}\subseteq A.
\]

\bigskip

In \S\ 2 below we review (Theorems K, K1, K2) the results we need, and prove
our main result, Theorem 1, extending Kendall's Theorem to the Karamata
setting, and Theorem 2, the Characterisation Theorem for `quasi regular
variation', where limits are taken avoiding a negligible exceptional set
(cf. BGT \S\ 2.9). In \S\ 3, we turn to `Beurling regular variation' [Ost2],
in its Baire version, proving the UCT in this setting (Th. 3) and the
relevant version of Kendall's Theorem (Th. 4), involving the functional
equation ($CJ$). We give the results we need on infinite combinatorics in
\S\ 4 (Prop. 3). In \S\ 5 we deal with \textit{general regular variation}
(Beurling setting, Baire versions). Here, the relevant version of Kendall's
Theorem (Th. 5) involves the functional equation ($BG$). Measure versions
(Th. 1M, Th. 4M, Th. 5M) follow in \S\ 6. Then in \S ~7 we turn to the
regular variation of the various sequences appearing in Kendall's Theorem
(Theorems 6I or 6M depending on context -- hereafter Theorem 6 for brevity).

Theorem 6). Character degradation (Theorem 7) resulting from ess-lim follows
in \S ~8 (cf. that from limsup and liminf in [BinO2]). Complements are
presented in \S\ 9 which we close with an Appendix on the relevant aspects
of coding (i.e. the links between classical and effective descriptive set
theory).

\bigskip

\textbf{2. Characterization theorems: the Karamata setting. }

Below $\mathbb{R}_{+}:=(0,\infty ),$ and functions are \textit{Baire} if
they have the Baire property, BP. We recall from [BinO1] that a divergent
sequence $c_{n}$ (i.e. with $\lim \sup_{n\rightarrow \infty }c_{n}=\infty $
) is said to be \textit{additively admissible}, resp. \textit{%
multiplicatively admissible, }if \textit{\ }%
\[
\lim \sup_{n\rightarrow \infty }c_{n+1}-c_{n}=0,\qquad \text{resp. }\lim
\sup_{n\rightarrow \infty }c_{n+1}/c_{n}=1.
\]%
As usual, we pass between multiplicative and additive versions at will by
using the exp/log isomorphisms between the additive group $\mathbb{R}$ (Haar
measure = Lebesgue measure $\mathrm{d}x)$ and the multiplicative group $%
\mathbb{R}_{+}$ (Haar measure $\mathrm{d}x/x);$ cf. [BGT, Ch. 1].

\bigskip

\noindent \textbf{Theorem K (Characterization theorem of Karamata regular
variation}, cf. [BGT, 1.4.1]). \textit{If }$f:\mathbb{R}_{+}\rightarrow
\mathbb{R}_{+}$\textit{\ is Baire/measurable and regularly varying, that is
for some function }$g$%
\[
\lim\nolimits_{x\rightarrow \infty }f(tx)/f(x)=g(t)\qquad (\forall t>0),
\]%
\textit{then }$g$\textit{\ is Baire/measurable and multiplicative:}%
\begin{equation}
g(st)=g(s)g(t)\qquad (\forall s,t>0),  \tag{$CFE$}
\end{equation}%
\textit{and so for some} $\gamma \in \mathbb{R}$%
\[
g(t)=t^{\gamma }.
\]

\noindent \textbf{Theorem K1 (Kingman's Croftian Theorem }[Kin1,2], cf.
[BGT, 1.9.1]).\textbf{\ }\textit{Take }$\{c_{n}\}$\textit{\ additively
admissible, }$I$\textit{\ an open interval of }$\mathbb{R}.$

\noindent (i)\textit{\ If }$G\subseteq \mathbb{R}$\textit{\ open and
unbounded from above, then}%
\[
c_{n}+x\in G\text{ infinitely often}
\]%
\textit{for some }$x\in I$\textit{.}

\noindent (ii)\textit{\ If }$f:\mathbb{R}\rightarrow \mathbb{R}$\textit{\
continuous} \textit{and}%
\[
\lim\nolimits_{n\rightarrow \infty }f(c_{n}+x)\text{ exists for all }x\in I,
\]%
\textit{then}%
\[
\lim\nolimits_{x\rightarrow \infty }f(x)\text{ exists.}
\]

\bigskip

The next result is Kendall's sequential characterization theorem of regular
variation.

\bigskip

\noindent \textbf{Theorem K2 (Kendall's Theorem }[Ken, Th. 16], cf. [BGT,
1.9.2]).\textbf{\ }\textit{For }$\{x_{n}\}_{n\in \mathbb{N}}$\textit{\
multiplicatively admissible and }$f:\mathbb{R}_{+}\rightarrow \mathbb{R}_{+}$%
\textit{\ continuous: if, as }$n\rightarrow \infty $\textit{,}%
\[
a_{n}f(\lambda x_{n})\rightarrow g(\lambda )\qquad (\lambda \in I)
\]%
\textit{for some interval }$I\subseteq (0,\infty ),$\textit{\ positive
sequence }$\{a_{n}\}_{n\in \mathbb{N}}$ \textit{and continuous function }$%
g:I\rightarrow \mathbb{R}_{+}$\textit{, then }$f$\textit{\ is regularly
varying: for each }$t>0,$%
\[
K(t):=\lim\nolimits_{x\rightarrow \infty }f(tx)/f(x)
\]%
\textit{exists, is finite, multiplicative, and both Baire and measurable. So
}$K(t)=t^{\kappa }$\textit{\ for some }$\kappa .$

\bigskip

The interval $I$ here may be arbitrarily small: a smidgen's worth of
sequential regular variation implies true regular variation of $f,$ and $%
\{a_{n}\}_{n\in \mathbb{N}}$:

\bigskip

\noindent \textbf{Corollary }(cf. Theorem 6I or 6M, \S ~7)\textbf{. }\textit{%
In Kendall's Theorem, }$\{a_{n}\}_{n\in \mathbb{N}}$\textit{\ is regularly
varying relative to }$\{x_{n}\}_{n\in \mathbb{N}}$\textit{\ with index }$%
-\kappa :$ \textit{if }$f(x)\sim x^{\kappa }\ell (x),$ \textit{with }$\ell $%
\textit{\ slowly varying, then for some constant }$c$%
\[
a_{n}\sim cx_{n}^{-\kappa }/\ell (x_{n}).
\]

Intervals as such are not needed here: the same is true for arbitrarily
small non-negligible (Baire/measurable) sets (see \S\ 7). Notice that, given
the hypotheses above, the limit function $K(t)$ is in fact the sequential
limit $\lim\nolimits_{n\rightarrow \infty }f(nt)/f(n)$ (of Baire/measurable
functions in the former case, and continuous functions in the latter), so is
Baire/measurable. The final assertions, characterizing the relevant limit
function, follow from theorems concerning Baire/measurable solutions of ($%
CFE $), the Cauchy functional equation (see e.g. [BinO3]).

Variants on the characterization theorem above are possible. First, one may
drop any condition of `topological good behaviour' (BP, or measurability,
which is BP under a change from the Euclidean to the density topology; see
e.g. [BinO6]), and weaken the quantifier on $t$ above, at the cost of
imposing a side-condition (the classical prototype is the Heiberg-Seneta
condition: BGT, 1.4.3) -- see [BinO8, \S\ 7]. By contrast, here we take the
passage to the limit \textit{sequentially }as in Kendall's Theorem\textit{, }%
through a suitable (admissible) sequence $\{x_{n}\}$, with our function $f$
again appearing \textit{once, }rather than \textit{twice}, but allow
exceptions on a meagre set. Our conclusion is of regular variation \textit{%
off an exceptional set -- `quasi regular variation', }as we shall call it.
This reduces to ordinary regular variation if we require also `topological
good behaviour' of the `essential limit' (below). As in [BinO2, BinO5, \S\ %
11], the passage to the essential limit results in \textit{character
degradation }(\S\ 7)\textit{. }Examination of this requires us to specify
the set-theoretic axioms we use (cf. [BinO2, 9]).

Our first definition covers both category and measure needs, again by
passage to the density topology.

\bigskip

\noindent \textbf{Definition 1. }For $L_{f}\ $finite, say that $f(x)$
\textit{has essential limit} $L_{f}$ as $x\rightarrow \infty $ and write $%
f(x)\rightarrow ^{\text{ess}}L_{f}$, or $\mathrm{ess}$-$\lim_{x\rightarrow
\infty }f(x)=L_{f},$ if for each $\varepsilon >0$ there is $X_{\varepsilon
}^{f}\in \mathbb{R}$ and meagre $M_{\varepsilon }^{f}$ such that%
\begin{equation}
|f(x)-L_{f}|<\varepsilon \text{ for all }x>X_{\varepsilon }^{f}\text{ off
the set }M_{\varepsilon }^{f}.  \tag{$\ast $}
\end{equation}

\noindent \textbf{Definition 2. }Say that a Baire function \textit{\ }$f:%
\mathbb{R}_{+}\rightarrow \mathbb{R}_{+}$ is \textit{quasi regularly varying}
weakly (resp. strongly) if $g(t):=\mathrm{ess}$-$\lim\nolimits_{x\rightarrow
\infty }f(tx)/f(x)$ exists and is finite (and resp. $g$ is Baire).

\bigskip

Our new results here are variants on or additions to Kendall's Theorem,
particularly Theorems 1 and 6 below.

\bigskip

\noindent \textbf{Theorem 1. }\textit{For }$\{x_{n}\}_{n\in \mathbb{N}}$%
\textit{\ multiplicatively admissible and }$f:\mathbb{R}_{+}\rightarrow
\mathbb{R}_{+}$\textit{\ Baire, if }%
\[
a_{n}f(\lambda x_{n})\rightarrow g(\lambda )\qquad (\lambda \in B)
\]%
\textit{for some non-meagre Baire set }$B\subseteq (0,\infty ),$\textit{\
positive sequence }$\{a_{n}\}_{n\in \mathbb{N}}$ \textit{and function }$%
g:B\rightarrow \mathbb{R}_{+}$\textit{, then }$f$\textit{\ is }(\textit{%
strongly})\textit{\ quasi regularly varying: for each }$s>0,$%
\[
K(s):=\mathrm{ess}\text{-}\lim\nolimits_{\lambda \rightarrow \infty
}f(s\lambda )/f(\lambda )
\]%
\textit{exists and is finite, and multiplicative. As }$g$ \textit{is Baire
on }$B$\textit{, }$K$ \textit{is locally bounded near }$s=1$\textit{, and so
}$K(s)=s^{\kappa }$\textit{\ for some }$\kappa \in \mathbb{R}.$

\bigskip

Some such result was suggested by [Bin2, footnote p. 162] in a discussion of
Kendall's Theorem. The question arises of strengthening Theorem 1 by
`thinning': requiring convergence for a smaller $\lambda $-set. Such
`quantifier weakening' is possible, and involves results of Steinhaus-Weil
type; see [BinO7,8] and \S\ 9.5. We delay the proof of Theorem 1 to
establish some preparatory results.

\bigskip

\noindent \textbf{Lemma 1. }(i)\textbf{\ }Essential limits preserve sums: if
$f(x)\rightarrow ^{\text{ess}}L_{f}$ and $g(x)\rightarrow ^{\text{ess}%
}L_{g}, $ then $(f+g)(x)\rightarrow ^{\text{ess}}L_{f}+L_{g};$ \textit{%
likewise for products.}

\noindent (ii) \textit{If }$h(x+u)-h(x)\rightarrow ^{\text{ess}}L_{u}$
\textit{and} $h(x+v)-h(x)\rightarrow ^{\text{ess}}L_{v},$ \textit{then}%
\[
h(x+u+v)-h(x)\rightarrow ^{\text{ess}}L_{u}+L_{v}.
\]

\bigskip

\noindent \textbf{Proof. }(i) For $\varepsilon >0$ choose $X_{\varepsilon
}^{f}\in \mathbb{R}$ and meagre $M_{\varepsilon }^{f}$ and likewise $%
X_{\varepsilon }^{g}\in \mathbb{R}$ and meagre $M_{\varepsilon }^{g}$ so
that ($\ast $) above holds for $f$ and $g$ respectively. Then ($\ast )$ for $%
(f+g)$ holds (with $2\varepsilon $ in lieu of $\varepsilon )$ for all $%
x>\max \{X_{\varepsilon }^{f},X_{\varepsilon }^{g}\}$ off the meagre set $%
M_{\varepsilon }^{f}\cup M_{\varepsilon }^{g}.$ Logarithmic transformation
yields the analogous result for products.

\noindent (ii) With $f(x):=h(x+u)-h(x)$ and $g(x):=h(x+v)-h(x),$ since $%
f(x)\rightarrow ^{\text{ess}}L_{u}$ and $g(x)\rightarrow ^{\text{ess}}L_{v}$,%
\[
\lbrack h(x+u+v)-h(x+v)]+[h(x+v)-h(x)]\rightarrow ^{\text{ess}}L_{u}+L_{v},
\]%
that is%
\[
\lbrack h(x+u+v)-h(x)]\rightarrow ^{\text{ess}}L_{u}+L_{v}.\qquad \square
\]

\bigskip

\noindent \textbf{Corollary. }\textit{For }$h:\mathbb{R}\rightarrow \mathbb{R%
}$ \textit{Baire and }$k:\mathbb{R}\rightarrow \mathbb{R}$\textit{\
arbitrary, }$\mathbb{G}_{\text{ess}}:=\{u:h(x+u)-h(x)\rightarrow ^{\text{ess}%
}k(u)\}$\textit{\ is a subgroup, and}%
\[
k(u+v)=k(u)+k(v)\qquad (u,v\in \mathbb{G}_{\text{ess}}).
\]%
\textit{So if }$\mathbb{G}_{\text{ess}}$\textit{\ contains a non-meagre
Baire set, then }$\mathbb{G}_{\text{ess}}=\mathbb{R}$\textit{\ and, if }$k$
\textit{is Baire, then }$k$\textit{\ is linear: }$k(u)=cu$\textit{.}

\bigskip

\noindent \textbf{Proof. }That $\mathbb{G}_{\text{ess}}=\mathbb{R}$ follows
here from the Subgroup Theorem for category [BGT, Cor. 1.1.4]; evidently $%
0\in \mathbb{G}_{\text{ess}},$ so it suffices to note that $-u\in \mathbb{G}%
_{\text{ess}}$ for $u\in \mathbb{G}_{\text{ess}}$: indeed with $y=x-u$%
\[
h(x-u)-h(x)=-[h(y+u)-h(y)]\rightarrow ^{\text{ess}}-k(u).\qquad \square
\]

\bigskip

\noindent \textbf{Theorem 2 } (\textbf{Characterization of quasi regular
variation).}\textit{\ If }$f:\mathbb{R}_{+}\rightarrow \mathbb{R}_{+}$%
\textit{\ is Baire/measurable and weakly quasi regularly varying with
essential limit function }$g,$%
\[
\mathrm{ess}\text{-}\lim\nolimits_{x\rightarrow \infty
}f(tx)/f(x)=g(t)\qquad (\forall t>0),
\]%
\textit{then }$g$\textit{\ is multiplicative. Furthermore, if }$g$ \textit{%
is Baire (i.e. }$f$\textit{\ is strongly quasi regularly varying), then for
some} $\gamma $%
\[
g(t)=t^{\gamma }.
\]

This follows from Lemma 1(ii) and its Corollary. As to the assumption that $%
g $ is Baire, Theorem 7 (in \S\ 8) clarifies the topological character of $g$%
. We turn now to a stronger form of Theorem K1(ii), which is based on our
generalizations of Theorem K1(i) in [BinO1].

\bigskip

\noindent \textbf{Proposition 2.} \textit{For }$f$ \textit{Baire and }$%
\{c_{n}\}_{n\in \mathbb{N}}$\textit{\ additively admissible,\ if }$%
\lim_{n}f(c_{n}+x)$\textit{\ exists for each }$x$\textit{\ in a Baire
non-meagre set }$C$\textit{, then }$\mathrm{ess}$-$\lim_{x\rightarrow \infty
}f(x)$\textit{\ exists, and, for quasi all }$x\in C,$\textit{\ equals }$%
\lim_{n}f(c_{n}+x)$\textit{. }

\bigskip

\noindent \textbf{Proof. }This follows [BGT, Th. 1.9.1(ii)]. W.l.o.g. $%
C:=I\backslash M$ with $I$ an open interval and $M$ meagre. The function $%
\hat{f}(x):=\lim_{n}f(c_{n}+x)$ is Baire on $C$. By the Baire-Kuratowski
Continuity Theorem (see e.g. [Oxt, Th. 8.1]) w.l.o.g. (expanding $M$ as
necessary) $\hat{f}|(I\backslash M)$ is continuous. Fix $x_{0}\in
I\backslash M;$ then, for any $\varepsilon >0,$ the set%
\[
J_{\varepsilon }:=\{x\in C:|\hat{f}(x)-\hat{f}(x_{0})|<\varepsilon \}
\]%
is open relative to $I\backslash M.$ So wlog $J_{\varepsilon }\subseteq I$
and $|\hat{f}(x)-\hat{f}(x_{0})|<\varepsilon $ holds on quasi all of $%
J_{\varepsilon }.$

We show that $f(x)\rightarrow ^{\text{ess}}\hat{f}(x_{0}).$ Otherwise, for
some $\varepsilon >0$ the Baire set%
\[
H:=\{x:|f(x)-\hat{f}(x_{0})|\geq \varepsilon \}
\]%
is essentially unbounded. Hence, by a generalization of Theorem K1 [BinO1,
Th. 3.6C], for quasi all $x\in J_{\varepsilon }$ there are infinitely many $%
n $ with $c_{n}+x\in H,$ i.e. $|f(c_{n}+x)-\hat{f}(x_{0})|\geq \varepsilon $
for infinitely many $n$ . For any such fixed $x\in J_{\varepsilon },$
passing to the limit yields $|\hat{f}(x)-\hat{f}(x_{0})|\geq \varepsilon ,$
and so this holds on quasi all of $J_{\varepsilon }.$ This contradicts the
reverse inequality, which holds on almost all of $J_{\varepsilon }$. $%
\square $

$\bigskip $

\noindent \textbf{Proof of Theorem 1. }We work in the multiplicative
positive reals, and begin by recalling a Kemperman-type Displacements Lemma
([BinO1, Cor. p. 157] -- there in additive notation) asserting that for $B$
Baire non-meagre, $B\cap sB$ is non-meagre for all $s$ close enough to $1$
-- for $s\in J_{\varepsilon }:=((1+\varepsilon )^{-1},1+\varepsilon ),$ say,
for some $\varepsilon >0$, cf. Theorem 6 below. (This may also be deduced
from the Pettis-Piccard Theorem, [Pet], [Pic], [BinGT, Th. 1.1.1], [BinO3,
Th. P].) Since scaling preserves category, this implies that $C(s):=B\cap
s^{-1}B$ is non-meagre for any $s\in J_{\varepsilon }.$

Now, for the most part, we follow Kendall's proof [Ken, Th. 16, p. 192]: if $%
\lambda \in B$ and $s\lambda \in B$ (i.e. $\lambda \in C(s)$), then%
\[
f(s\lambda x_{n})/f(\lambda x_{n})=a_{n}f(s\lambda x_{n})/a_{n}f(\lambda
x_{n})\rightarrow g(s\lambda )/g(\lambda ).
\]%
Put $k_{s}(\lambda ):=f(s\lambda )/f(\lambda );$ then for any fixed $s\in
J_{\varepsilon }$, for $\lambda $ in the Baire non-meagre set $C(s),$
\[
k_{s}(\lambda x_{n})\rightarrow g(s\lambda )/g(\lambda ),\text{ as }%
n\rightarrow \infty .
\]%
By Prop. 2, $K(s):=\mathrm{ess}$-$\lim_{x\rightarrow \infty }k_{s}(x)$
exists for this arbitrary $s\in J_{\varepsilon };$ furthermore, for each
such $s$ and quasi all $\lambda \in C(s),$%
\[
K(s)=\lim\nolimits_{n}f(s\lambda x_{n})/f(\lambda x_{n})=g(s\lambda
)/g(\lambda ).
\]

As in the Corollary above, $\mathbb{G}:=\{s:K(s)$ is well defined$\}$ is a
multiplicative subgroup, since, for fixed $s,t\in \mathbb{G}$,%
\begin{eqnarray*}
K(st) &=&\mathrm{ess}\text{-}\lim_{x\rightarrow \infty }k_{st}(x)=\mathrm{ess%
}\text{-}\lim_{x\rightarrow \infty }[f(stx)/f(tx)]\cdot \mathrm{ess}\text{-}%
\lim_{x\rightarrow \infty }[f(tx)/f(x)] \\
&=&K(s)K(t).
\end{eqnarray*}%
Moreover, $\mathbb{G}$ contains the interval $J_{\varepsilon },$ so, by the
Steinhaus Subgroup Theorem (see e.g. [BinO2, Th. 6.2]), $\mathbb{G=R}_{+}$.
So $\mathrm{ess}$-$\lim_{x\rightarrow \infty }k_{s}(x)$ exists for all $s$
and is multiplicative, as in Lemma 1(ii).

For the last part, being a sequential limit $g$ is (positive and) Baire. By
passing to a smaller non-meagre subset of $B,\ $we may w.l.o.g. assume that $%
g$ is bounded on $B,$ so that, for some $0<a<b,$ say:%
\[
a<g(\lambda )<b\qquad (\lambda \in B).
\]%
Likewise, by passage to a corresponding smaller $\varepsilon >0,$ if
necessary, we again conclude (as above) that, for $s\in J_{\varepsilon }$
and quasi all $\lambda \in C(s),$%
\[
K(s)=g(s\lambda )/g(\lambda )>0.
\]%
For $s\in J_{\varepsilon }$ choose $\lambda _{s}\in C(s)$ to witness the
preceeding equation. Then%
\[
K(s)=g(s\lambda _{s})/g(\lambda _{s})\in (a/b,b/a).
\]%
So $K$ is locally bounded on $J_{\varepsilon }.$ Hence, by the Darboux
Theorem $($[Dar1,2], [Kuc3, \S 14.4]) or the Banach-Mehdi Theorem (see e.g.
[BinO3, Th. BM]), $K$ is continuous, and so a power function: $%
K(s)=s^{\kappa }$ for all $s>0,$ since $K>0.$ $\square $

\bigskip

\noindent \textbf{Remarks. 1 (}\textit{`Nice' versions of Baire functions}%
\textbf{). }The proof above centers on $k(s,t):=g(st)/g(t)$ as a function of
two variables, with $g$ Baire. It is instructive to take note why the
composite function $k$ may be assumed Baire.

In order for $f:\mathbb{R}^{d}\mathbb{\rightarrow R}$ to be a Baire function
it is necessary and sufficient that the restriction $f|(\mathbb{R}\backslash
M)$ be continuous for some meagre set $M$ [Oxt, Th. 8.1]. Above we opted to
neglect behaviour on meagre sets, so $M$ may as well be expanded to a union
of closed nwd (nowhere dense) sets, and so for $f$ Baire, the restriction $%
f|H_{f}$ is continuous for some dense $\mathcal{G}_{\delta }$ set $H_{f}.$
With this in mind, note why the map $(s,t)\mapsto f(st)$ may be taken Baire:
for each $q\in \mathbb{Q}$, the level set $L(q):=\{(s,t):f(st)<q\}$ is the
projection of the set%
\[
\{(s,t,u):u=st\}\cap \{(s,t,u):(s,t)\in H_{f}\times H_{f}\text{ \& }%
f(u)<q\};
\]%
as the second term here is a $\mathcal{G}_{\delta },$ the projection $L(q)$
is analytic, and so by Nikodym's Theorem [Rog, Cor. 2.9.4] [Kec, 29.14] is
again Baire. As $\mathbb{Q}$ is countable, all the rational level sets are $%
\mathcal{G}_{\delta }$ sets modulo one single union of closed nwd (nowhere
dense) sets.

The moral is: we may pass to a `version' of $k$ which is `nice': all its
rational level sets are $\mathcal{G}_{\delta }$. Of course, the actual
function $k$ can be as `nasty' as the meagre set one's set-theory admits,
and that depends on one's selection of a (perhaps weak) form of the axiom of
choice. See e.g. the Appendix below and [BinO9].

\noindent \textbf{2. }In the last part of our proof, where we deduce the
form of $K$, we necessarily parted company with Kendall's proof which at
that point refers to Theorem K to identify $K$. For, Theorem 1 refers to
essential, rather than ordinary, limits; nor could we apply the Corollary
above to $k=\log K$, as we did not then know the topological character of $%
K. $

\bigskip

\textbf{3. Beurling regular variation: Baire versions.}

\textit{Beurling slow variation} [BGT, \S\ 2.11], relative to a function $%
\varphi ,$ was used by Beurling to prove a Tauberian theorem for Borel
summability, which is not of convolution form but `convolution-like' (see
e.g. [Bin1]; for links with Riesz means, see [Bin3]). For \textit{Beurling's
Tauberian theorem, } extending Wiener's Tauberian Theorem, see e.g. [Kor,
IV.11], [BinO10, \S\ 6.1] (there, the `convolution-like' operation is shown
to be an `asymptotic convolution'). This has led recently to a
generalization [BinO5,10] of RV\ to $\varphi $-RV, shortly to be recalled.
Below, Kendall's Theorem is first extended to this context. It has emerged
recently that the most convenient way to define Beurling's idea is to use
some simple algebraic tools. Key here is notation introduced by Popa [Pop]
(and later independently by Javor [Jav]) to study the equation ($GS$) above,
whose central role for RV was established only quite recently. For arbitrary
$h:\mathbb{R}\rightarrow \mathbb{R}$ define $\circ _{h},$ the (\textit{Popa}%
) \textit{circle operation }[BinO6] on $\mathbb{R}$\textit{, }by%
\[
s\circ _{h}t=s+th(s).
\]%
When $h(t)=\eta _{1}(t)\equiv 1+t$ this reduces to the circle operation $%
s\circ t=s+t+st$ of ring theory (for background see [BinO6, \S\ 3], [Ost3,
\S\ 2]). In the case of $\mathbb{R}$ the operation endows it with a group
structure conjugate to ordinary multiplication in view of the identity%
\[
s\circ t=s+t+st=(1+s)(1+t)-1.
\]%
For $\eta $ satisfying $GS,$ we denote the resulting \textit{Popa (circle)
group} by $\mathbb{G}_{\eta }^{\ast }:=\{x:\eta (x)\neq 0\},$ writing its
inverse operation as $x_{\eta }^{-1}.$

In this notation, as above, a function $\varphi :\mathbb{R}_{+}\rightarrow
\mathbb{R}_{+}$ is said to be \textit{self-equivarying}, $\varphi \in SE,$
if it is $O(x)$ as $x\rightarrow \infty $ and
\[
\eta _{x}(t):=\varphi (x\circ _{\varphi }t)/\varphi (x)\rightarrow \eta
(t)\qquad (\text{locally uniformly in }t)\text{.}
\]%
By Theorem O the limit function satisfies the \textit{Go\l \k{a}b-Schinzel
equation} of \S ~1. So $\circ _{\eta }$ is commutative and associative. We
recall that positive solutions, relevant here, are of the form $\eta
(t)=\eta _{\rho }(t):=1+\rho t,$ for $t>-1/\rho $ with $\rho \geq 0,$ and
that the case $\rho =0$ when $\eta \equiv 1$ with $\varphi $ an $o(x)$
function corresponds to Beurling's original notion of a \textit{%
self-neglecting} function.

\bigskip

\noindent \textbf{Definition 3. }For $\varphi \in SE$ a function $f:\mathbb{R%
}_{+}\rightarrow \mathbb{R}_{+}$ is $\varphi $\textit{-regularly varying}
if, for some $g$ and all $t>0,$
\[
f(x+t\varphi (x))/f(x)\rightarrow K(t)\qquad (\text{locally uniformly in }t)%
\text{;}
\]%
that is, $f(x\circ _{\varphi }t)/f(x)\rightarrow K(t).$

\bigskip

Here again, as in the Karamata setting of regular variation of \S\ 1 the%
\textit{\ Uniform Convergence Theorem} holds [Ost2, Th. 1]: if $\varphi \in
SE,f,K$ are all Baire/measurable, then convergence is necessarily locally
uniform. Our next result extends the UCT to our present Kendall setting, and
justifies below the assumption of local uniformity throughout. Here $g,$ as
a sequential limit is Baire/measurable.

\bigskip

\noindent \textbf{Theorem 3 (UCT}, cf. [Ost2, Th. 1] and [BinO5, Th. 2B/M]).
\textit{Take }$B$\textit{\ non-negligible, }$\varphi \in SE,$\textit{\ }$f$ (%
\textit{and so} $g:B\rightarrow \mathbb{R}_{+}$) \textit{all Baire or all
measurable. If }%
\[
a_{n}f(x_{n}\circ _{\varphi }t)\rightarrow g(t)\qquad (t\in B),
\]%
\textit{then}
\[
a_{n}f(x_{n}\circ _{\varphi }t)\rightarrow g(t)\qquad (t\in B)\qquad (\text{%
locally uniformly in }t).
\]

\noindent \textbf{Proof.} This follows from Prop. 1 (Affine Two-sets Lemma,
above) as in the proof of [Ost2, Th. 1] (cf. [BinO5, Th. 2B/M]) with the
following changes. First, replace $h_{N}(x\circ _{\varphi }t)$ by $\log
f(x\circ _{\varphi }t)$ (`N for numerator'). Next, replace $h_{D}(x_{n})$
(`D for denominator') by $\log (1/a_{n}).$ Finally, replace $\mathbb{R}$ by $%
B.$ $\square $

\bigskip

\noindent \textbf{Blanket Assumption of non-triviality on }$g$\textbf{. }We
will call a function $g:B\rightarrow \mathbb{R}_{+}$, as above, \textit{%
trivial} if it takes values only in $\{0,1\}.$

\bigskip

In the corresponding functional equations in $K\ $that arise below, this in
turn excludes trivial solutions (note that the multiplicative Cauchy
functional equation has trivial solutions identically $0$ or $1$, the first
excluded in Theorem 1): see \S\ 9.6.

\bigskip

\noindent \textbf{Theorem 4.} \textit{Take }$B$\ \textit{Baire non-meagre, }$%
\{x_{n}\}_{n\in \mathbb{N}}$\textit{\ additively admissible, and }$\varphi
\in SE$\textit{\ and} $f$ \textit{both} \textit{Baire. If the Kendall
condition }%
\[
a_{n}f(x_{n}\circ _{\varphi }t)\rightarrow g(t)\qquad (t\in B)\qquad (\text{%
locally uniformly in }t)
\]%
\textit{holds, then }$f$\textit{\ is }$\varphi $\textit{-regularly varying
and}%
\[
f(x+t\varphi (x))/f(x)\rightarrow K(t),\text{ for }t>0,
\]%
\textit{with }$K$\textit{\ continuous, satisfying\ the Chudziak-Jab\l o\'{n}%
ska equation:}
\begin{equation}
K(u\circ _{\eta }v)=K(v)K(u).  \tag{$CJ$}
\end{equation}

\bigskip

\noindent \textbf{Proof.} Put $h_{n}(s):=(x_{n}\circ _{\varphi }\lambda
)\circ _{\varphi }s$ and note that%
\begin{eqnarray*}
h_{n}(s) &=&(x_{n}+\lambda \varphi (x_{n}))+s\varphi (x_{n}+\lambda \varphi
(x_{n})) \\
&=&x_{n}+\varphi (x_{n})[\lambda +s\eta _{x_{n}}(\lambda )]=x_{n}\circ
_{\varphi }[\lambda +s\eta _{x_{n}}(\lambda )].
\end{eqnarray*}%
Also, for fixed $\lambda $ and $n,$ the function $s\mapsto h_{n}(s)$ is a
homeomorphism of $\mathbb{R}$ under the usual (Euclidean) topology; as $%
\varphi (x)$ is $O(x),$%
\[
(x_{n+1}\circ _{\varphi }\lambda )-(x_{n}\circ _{\varphi }\lambda
)=(x_{n+1}-x_{n})+\lambda (\varphi (x_{n+1})-\varphi (x_{n}))\rightarrow 0,
\]%
and so -- since $\varphi (x_{n}+\lambda \varphi (x_{n}))/\varphi
(x_{n})\rightarrow \eta (\lambda )$ -- the (bitopological) generalization of
Theorem K1 [BinO1, Th. 3.5] holds for the homeomorphisms $\{h_{n}\}_{n\in
\mathbb{N}}.$

Now put%
\[
k_{s}(t):=f(t\circ _{\varphi }s)/f(t).
\]%
Then by the local uniformity in the Kendall condition of the theorem, because%
\[
\lambda +s\eta _{x_{n}}(\lambda )\rightarrow \lambda +s\eta (\lambda
)=\lambda \circ _{\eta }s,
\]%
with $\eta $ Baire (so continuous), we get
\[
a_{n}f(x_{n}\circ _{\varphi }[\lambda +s\eta _{x_{n}}(\lambda )])\rightarrow
g(\lambda \circ _{\eta }s).
\]%
As before, but now working in $\mathbb{G}_{\eta },$ the Popa group under $%
\circ _{\eta },$ there is an interval $J_{\varepsilon }$ of values $s$ for
which $C(s):=B\cap (s_{\eta }^{-1}\circ _{\eta }B)$ is non-meagre, and so
for quasi all $\lambda \in C(s)$%
\begin{eqnarray*}
k_{s}(x_{n}\circ _{\varphi }\lambda ) &=&f((x_{n}\circ _{\varphi }\lambda
)\circ _{\varphi }s)/f(x_{n}\circ _{\varphi }\lambda ) \\
&=&f(x_{n}+\varphi (x_{n})[\lambda +s\eta _{x}(\lambda )])/f(x_{n}\circ
_{\varphi }\lambda ) \\
&\rightarrow &g(\lambda \circ _{\eta }s)/g(\lambda ).
\end{eqnarray*}%
Then, as before, by Prop. 2, $K(s):=\mathrm{ess}$-$\lim_{x\rightarrow \infty
}k_{s}(x)$ exists for this arbitrary $s\in J_{\varepsilon }.$ Also, for each
such $s$ and quasi all $\lambda \in C(s),$%
\[
K(s)=\lim\nolimits_{n}k_{s}(x_{n}\circ _{\varphi }\lambda )=g(\lambda \circ
_{\eta }s)/g(\lambda )>0.
\]%
and
\[
K(s)=\lim\nolimits_{n}f((x_{n}\circ _{\varphi }\lambda )\circ _{\varphi
}s)/f(x_{n}\circ _{\varphi }\lambda )=g(\lambda \circ _{\eta }s)/g(\lambda
).
\]%
Since $\eta =\eta _{\rho },$ for some $\rho \geq 0,$ by the Steinhaus
subgroup theorem, as in the Corollary, $K(s)$ exists for all $s$ in the Popa
group $\mathbb{G}_{\eta }:=\{s:\eta (s)=1+\rho s\neq 0\},$ since $K\ $is a
homomorphism. Indeed, since%
\[
\tilde{v}(u,x):=v[\eta (u)/\eta _{x}(u)]\rightarrow v,
\]%
and%
\[
x\circ _{\varphi }(u\circ _{\eta }v)=(x\circ _{\eta }u)\circ _{\varphi }%
\tilde{v}(u,x)=(x+\varphi (x)u)+\tilde{v}(u,x)\varphi (x+\varphi (x)u)),
\]%
\begin{eqnarray*}
K(u\circ _{\eta }v) &=&\mathrm{ess}\text{-}\lim_{x\rightarrow \infty
}k_{u\circ _{\eta }v}(x)=\mathrm{ess}\text{-}\lim_{x\rightarrow \infty
}f(x\circ _{\varphi }(u\circ _{\eta }v))/f(x) \\
&=&\mathrm{ess}\text{-}\lim_{x\rightarrow \infty }[f((x\circ _{\eta }u)\circ
_{\varphi }\tilde{v}(u,x))/f(x\circ _{\eta }u)]\cdot \mathrm{ess}\text{-}%
\lim_{x\rightarrow \infty }[f(x\circ _{\eta }u)/f(x)] \\
&=&K(v)K(u),
\end{eqnarray*}%
giving ($CJ$). $\square $

\bigskip

\textbf{4. Croftian Infinite Combinatorics}

A key ingredient in the multiplicative form of Kingman's Theorem K1 is that,
for $\{d_{n}\}_{n\in \mathbb{N}}$ multiplicatively admissible, the sequence
of dilations (homeomorphisms) $h_{d}:x\mapsto dx$ for $d=d_{1},d_{2},...$
has the property that, for any non-degenerate interval $J=(a,b),$ the `tail
union' $\tbigcup\nolimits_{n\geq m}h_{d_{n}}(J)$ contains an infinite
half-line (cf. [BinO1, Th. 3.2]). This property no longer holds when $J$ is
replaced by a non-null (closed) set, as in the example due to Roy Davies
[BinO1, Th. 4.6]. To circumvent this, one may replace the indexing set $%
\mathbb{N}$ with its natural order by the (countable) set $\mathbb{Q}_{+}$
of positive rationals with their natural order (induced from the reals) and
employ the corresponding rational dilations $h_{q}(x)=qx.$ Then, for any
density-open set $W$ (see below), the corresponding tail union $%
\tbigcup\nolimits_{q\geq r}h_{q}(W)$ contains almost all of an infinite
half-line [BinO1, Th. 3.2]. Further [BinO1, Th. 3.2, Remark 2], this
continues to hold with $\mathbb{Q}_{+}$ replaced by any set of dilations $%
\{h_{d}(x):d\in D\}$ with $D$ dense in $\mathbb{R}_{+}$ (equivalently:
translations on $\mathbb{R}$). This may be read as saying that for fixed $x,$
the set $\{h_{d}(x):d\in D\}$ is dense in $\mathbb{R}_{+}$. We will also
need a strengthening of this provided in the next result, which uses the $%
\mathbb{Q}_{+}$ analogue of the homeomorphisms of Th. 4 above: for fixed
positive $q$ and $\lambda ,$ put%
\[
h_{q}(s):=(q\circ _{\varphi }\lambda )\circ _{\varphi }s=q+\lambda \varphi
(q)+s\varphi (q+\lambda \varphi (q)),
\]%
to be called the `$\varphi $-\textit{dilations}'; below we assume $\varphi $
is continuous and again rely on Prop. 1 (but cf. \S\ 9.7).

\bigskip

We recall that a set is \textit{density-open} if all its points are
(Lebesgue) density points. To maintain category-measure duality, for
convenience we are content to adopt the following

\bigskip

\noindent \textbf{Definition 4. }Call a set \textit{category-open} if it
takes the form of an open set less a meagre set.

\bigskip

Whilst the correct duality stance would be to work bitopologically and use
the category and measure versions of Hashimoto topologies [BinO7], all we
need below is that the intersection of two sets of one of these two types is
again of the same type.

\bigskip

\noindent \textbf{Proposition 3.}\textit{\ For }$A,B$\textit{\
category/density open with }$B$\ \textit{unbounded and }$\varphi $ \textit{%
continuous with }$\varphi =O(x)$\textit{, there are arbitrarily large
rationals }$q$\textit{\ and points }$a_{q}\in A,$\textit{\ }$b_{q}\in B$%
\textit{\ with }%
\[
h_{q}(a_{q})=b_{q}.
\]%
\textit{Hence the tail union }$\tbigcup\nolimits_{q\geq r}h_{q}(A)$\textit{\
contains quasi all of an infinite half-line.}

\bigskip

\noindent \textbf{Proof.} We consider only the density-open case, as the
category-open case is similar but simpler. Our aim is to establish the
relation
\[
q+\lambda \varphi (q)+a_{q}\varphi (q+\lambda \varphi (q))=b_{q},
\]%
as above, or equivalently%
\[
1+a_{q}m_{\lambda }(q)=\frac{b_{q}}{q+\lambda \varphi (q)}\text{ for }%
m_{\lambda }(x):=\frac{\varphi (x+\lambda \varphi (x))}{x+\lambda \varphi (x)%
}>0.
\]%
As $\varphi (x)=O(x)$, $m_{\lambda }(x)$ remains bounded as $x\rightarrow
\infty ,$ say by $M(\lambda ).$ Fix $a\in A$ arbitrarily, then choose $b\in
B $ as large as desired with $b>1+aM(\lambda ).$ Since $\varphi (x)$ and $%
m_{\lambda }(x)$ are continuous in $x,$ and $x+\lambda \varphi
(x)=x(1+\lambda \varphi (x)/x)\rightarrow \infty ,$ there exists $x$ (which
is as large as desired for $b$ large enough) with%
\[
\frac{b}{x+\lambda \varphi (x)}=1+am_{x}(\lambda ):\qquad a=\frac{1}{%
m_{x}(\lambda )(x+\lambda \varphi (x))}(b-[x+\lambda \varphi (x)]).
\]%
Fix such an $x,$ and choose a rational sequence $q_{n}\rightarrow x.$ Again,
by the continuity of both $\varphi $ and $m_{\lambda }$,
\[
c_{n}:=\frac{1}{m_{q_{n}}(\lambda )[q_{n}+\lambda \varphi (q_{n})]}%
\rightarrow c:=\frac{1}{m_{x}(\lambda )[x+\lambda \varphi (x)]}>0.
\]%
So
\[
a/c=\left( b-\frac{1}{cm_{x}(\lambda )}\right) :\qquad a/c\in
B_{0}:=(A/c)\cap \left( B-\frac{1}{cm_{x}(\lambda )}\right) .
\]%
Here $B_{0}$ is dense-open, since $a/c$ is a density point both of $A/c$ and
of the translate of $B;$ furthermore, $cB_{0}\subseteq A.$ Put
\[
z_{n}:=c_{n}\left( \frac{1}{cm_{x}(\lambda )}-[q_{n}+\lambda \varphi
(q_{n})]\right) .
\]%
Then, since $cm_{x}(\lambda )=1/(x+\lambda \varphi (x))$ and $%
m_{q_{n}}(\lambda )[q_{n}+\lambda \varphi (q_{n})]=\varphi (q_{n}+\lambda
\varphi (q_{n})),$
\begin{eqnarray*}
z_{n} &=&\frac{1/m_{q_{n}}(\lambda )}{q_{n}+\lambda \varphi (q_{n})}\left(
[x+\lambda \varphi (x)]-[q_{n}+\lambda \varphi (q_{n})]\right) \\
&=&\frac{1}{\varphi (q_{n}+\lambda \varphi (q_{n}))}\left( [x-q_{n}]+\lambda
\lbrack \varphi (x)-\varphi (q_{n})]\right) \rightarrow 0.
\end{eqnarray*}

By the Affine Two Sets Lemma above (applied to $A$ and $B_{0}$, rather than $%
A$ and $B),$ for almost all $b_{0}^{\prime }\in B_{0}$ the sequence%
\[
c_{n}b_{0}^{\prime }+z_{n}\in A\text{ i.o.}
\]%
In particular, as $B_{0}\subseteq B-1/cm_{x}(\lambda ),$ there are $%
a^{\prime }\in A,$ $b_{0}^{\prime }\in B_{0},$ $b^{\prime }\in B$ with $%
b_{0}^{\prime }=[b^{\prime }-1/cm_{x}(\lambda )],$ and some $n$ with%
\[
c_{n}(b^{\prime }-\frac{1}{cm_{x}(\lambda )})+z_{n}=a^{\prime }.
\]%
Substituting for $z_{n}$ gives%
\[
a^{\prime }=c_{n}\left( b^{\prime }-\frac{1}{cm_{x}(\lambda )}\right)
+c_{n}\left( \frac{1}{cm_{x}(\lambda )}-[q_{n}+\lambda \varphi
(q_{n})]\right) ,
\]%
that is
\[
a^{\prime }=\frac{1/m_{q_{n}}(\lambda )}{q_{n}+\lambda \varphi (q_{n})}%
\left( b^{\prime }-[q_{n}+\lambda \varphi (q_{n})]\right) :\qquad
1+a^{\prime }m_{q_{n}}(\lambda )=\frac{b^{\prime }}{q_{n}+\lambda \varphi
(q_{n})}.
\]

The final assertion follows verbatim as in [BinO1, Th. 3.2]. $\square $

\bigskip

\textbf{5. General regular variation: Beurling-Baire versions.}

An analysis similar to that in Theorem 4 may be performed for the general
setting of regular variation of \S\ 1, i.e. with asymptotics defined by%
\[
\lbrack f(x+t\varphi (x))-f(x)]/h(x)\rightarrow g(t)\qquad (\text{locally
uniformly in }t)\text{.}
\]%
The third function here, $h$ (assumed positive), must satisfy (see the
Remark below)%
\[
h(x\circ _{\varphi }t)/h(x)\rightarrow r(t)\qquad (\text{locally uniformly
in }t)\text{,}
\]%
so that, by Theorem 4, $r(t)$ satisfies (CJ).

\bigskip

\noindent \textbf{Theorem 5.} \textit{Take }$B\subseteq (0,\infty )$\
\textit{Baire non-meagre, }$\{x_{n}\}_{n\in \mathbb{N}}$\textit{\ additively
admissible, and }$\varphi \in SE,$ $f$ (\textit{and so }$g:B\rightarrow
\mathbb{R}_{+}$) \textit{all Baire. If the Kendall condition }%
\[
a_{n}[f(x_{n}+t\varphi (x_{n})-f(x_{n}))/h(x_{n})\rightarrow g(t)\qquad (%
\text{locally uniformly in }t)\text{.}
\]%
\textit{holds, then}%
\[
\lbrack f(x+t\varphi (x))-f(x)]/h(x)\rightarrow K(t)\qquad (\text{locally
uniformly in }t)\text{,}
\]%
\textit{with }$K$ \textit{continuous, satisfying\ the Beurling-Goldie
equation}
\begin{equation}
K(u\circ _{\eta }v)=K(u)\circ _{\sigma }K(v),  \tag{$BG$}
\end{equation}%
\textit{and with the }$\sigma $\textit{\ in the }$\circ _{\sigma }$\textit{\
above satisfying}
\[
\sigma (K(u))=r(u):\qquad \sigma (t)=r(K^{-1}(t)).
\]

\noindent \textbf{Proof. }Here one takes%
\[
k_{s}(t):=[f(t\circ _{\varphi }s)-f(t)]/h(t).
\]%
The analysis is unchanged (with $\tilde{v}(u,x)\rightarrow v$ in the same
notation), but the `kernel function' $K$ now satisfies%
\begin{eqnarray*}
K(u\circ _{\eta }v) &=&\mathrm{ess}\text{-}\lim_{x\rightarrow \infty
}k_{u\circ _{\eta }v}(x)=\mathrm{ess}\text{-}\lim_{x\rightarrow \infty
}[f(x\circ _{\varphi }(u\circ _{\eta }v))-f(x)]/h(x) \\
&=&\mathrm{ess}\text{-}\lim_{x\rightarrow \infty }\{[f((x\circ _{\eta
}u)\circ _{\eta }\tilde{v}(u,x))-f(x\circ _{\eta }u)]/h(x\circ _{\eta
}u)\}\cdot \mathrm{ess}\text{-}\lim_{x\rightarrow \infty }\frac{h(x\circ
_{\eta }u)}{h(x)} \\
&&+\mathrm{ess}\text{-}\lim_{x\rightarrow \infty }[f(x\circ _{\eta
}u)-f(x)]/h(x) \\
&=&K(v)r(u)+K(u).
\end{eqnarray*}%
Here $K$ is a homomorphism between the two Popa groups $\mathbb{G}_{\eta }$
and $\mathbb{G}_{\sigma }$. $\square $

\bigskip

\noindent \textbf{Remark. }Notice that the preceeding equation holds \textit{%
if and only if} $h$ has the asymptotic behaviour specified above. We also
note that a non-zero $K$ will necessarily be monotone -- see [Ost3].

\bigskip

\textbf{6. Measure versions.}

The aim is to read off measure analogues of Theorems K2, 1, 4 and 5 by
replacing the Euclidean topology by the density topology. Mutatis mutandis,
the `density open' version of Prop. 3 allows precisely this.

\bigskip

\noindent \textbf{Proposition 2M (}cf. [BinO1, Th. 4.1]). \textit{For }$f$
\textit{measurable,\ if }$\lim_{q\in \mathbb{Q},q\rightarrow \infty }f(qx)$%
\textit{\ exists for each }$x$\textit{\ in a non-null measurable set }$B$%
\textit{, then }$\mathrm{ess}$-$\lim_{x\rightarrow \infty }f(x)$\textit{\
exists, and, for almost all }$x\in B,$\textit{\ equals }$\lim_{q\in \mathbb{Q%
},q\rightarrow \infty }f(qx)$\textit{. }

\bigskip

\noindent \textbf{Proof. }The same proof as in Prop. 2 works with $\hat{f}%
(x)=\lim_{q\in \mathbb{Q},q\rightarrow \infty }f(qx),$ which is measurable
on $B$. By the Luzin Continuity Theorem (see e.g. [Oxt, Th. 8.2]) we may
assume that $\hat{f}|B$ is continuous (otherwise pass to a non-null subset $%
B^{\prime }$ of $B$ on which this holds, removing a part of measure as small
as desired). From here the proof is the same save that `for infinitely many $%
n\in \mathbb{N}$' is replaced by `on some unbounded sequence of $q\in
\mathbb{Q}_{+}$'. $\square $

\bigskip

This allows for the following measure versions. The first needs only
ordinary dilations $\{h_{q}:q\in \mathbb{Q}_{+}\};$ the second needs the $%
\varphi $-dilations of \S\ 4.

\bigskip

\noindent \textbf{Theorem 1M. }\textit{For }$f:\mathbb{R}_{+}\rightarrow
\mathbb{R}_{+}$\textit{\ measurable, if }%
\[
a_{q}f(q\lambda )\rightarrow g(\lambda )\qquad (\lambda \in B)\qquad
(q\rightarrow \infty \text{ through }\mathbb{Q}_{+})
\]%
\textit{for some non-null measurable set }$B\subseteq (0,\infty )$ \textit{%
and function }$g:B\rightarrow \mathbb{R}_{+}$\textit{, then }$f$\textit{\ is
weakly almost regularly varying: for each }$s>0,$%
\[
K(s):=\mathrm{ess}\text{-}\lim\nolimits_{\lambda \rightarrow \infty
}f(s\lambda )/f(\lambda )
\]%
\textit{exists and is finite, and multiplicative. As }$g$ \textit{is
measurable on }$B$\textit{, }$K$ \textit{is locally bounded near }$s=1$%
\textit{, and so }$K(s)=s^{\kappa }$\textit{\ for some }$\kappa .$

\bigskip

\bigskip

\noindent \textbf{Theorem 4M.} \textit{Take }$\varphi \in SE$\textit{\
continuous, }$B\subseteq (0,\infty )$\textit{\ measurable non-null and }$f$ (%
\textit{and so} $g:B\rightarrow \mathbb{R}_{+}$) \textit{measurable. If the
Kendall condition }%
\[
a_{q}f(q\circ _{\varphi }t)\rightarrow g(t)\qquad (t\in B,\text{ }%
q\rightarrow \infty ,\text{ }q\in \mathbb{Q})\qquad (\text{locally uniformly
in }t)
\]%
\textit{holds, then }$f$\textit{\ is }$\varphi $\textit{-regularly varying
and}%
\[
f(x+t\varphi (x))/f(x)\rightarrow K(t),\text{ for }t>0,
\]%
\textit{with }$K$\ \textit{continuous, satisfying\ the equation:}
\begin{equation}
K(u\circ _{\eta }v)=K(v)K(u).  \tag{$CJ$}
\end{equation}

\noindent \textbf{Proof.} Put $h_{q}(s):=(q\circ _{\varphi }\lambda )\circ
_{\varphi }s$ and apply the density-open variant of Prop. 3. So the
(bitopological) generalization of Theorem K1 [BinO1, Th. 3.5M] holds for the
homeomorphisms $\{h_{q}\}_{q\in \mathbb{Q}_{+}}.$ $\square $

\bigskip

The measure analogue Theorem 5M of Theorem 5 follows similarly.

\bigskip

\textbf{7. Regular variation: sequences.}

Theorem 6I below (`I for interval') brings out a property of regular
variation in the sequence $\{a_{n}\}_{n\in \mathbb{N}}$ of Kendall's
Theorem. This (which is actually left implicit in [Ken]) is important in
various contexts, particularly in probability theory; see e.g. [Bin2] and
\S\ 9.3 and 9.4 below. The reduction below to a power function brings the
Kendall condition into alignment with one due to Seneta, cf. BGT \S ~1.9.3.
Below, without loss of generality we assume that $g$ is continuous on $B$
(by the Baire-Kuratowski or the Luzin Continuity theorems -- passing to a
smaller set, as necessary).

We first need to isolate from the proof of Kendall's Theorem all the
available information about the $g$ function ocurring in that theorem. As
usual `negligible' below refers to meagre/null sets.

\bigskip

\noindent \textbf{Definition. }For $B$ non-negligible, say that the
continuous function $g:B\rightarrow \mathbb{R}_{+}$ satisfies the\textit{\
Restricted Cauchy functional equation, }(\textit{Res}$CFE$) on $B,$ if%
\begin{equation}
g(s\lambda )=s^{\kappa }g(\lambda )\qquad (\forall \lambda \in (B\cap
Bs^{-1})\backslash M(s)\text{ \textit{with }}M(s)\in \mathcal{N}\text{%
\textit{, }}\forall s\text{\textit{\ with }}B\cap Bs^{-1}\notin \mathcal{N}).
\tag{\QTR{it}{Res}$CFE_{B}$}
\end{equation}

This is novel here; for (conditional) functional equations, i.e. with
restricted domains see [AczD, Ch. 6, 7, 16], [Kuc2], [Kuc3, \S 13.6], and
for further background literature [Ree].

\bigskip

\noindent \textbf{Theorem 6I.} (i) \textit{For }$I$ \textit{an interval,} $M$
\textit{negligible and }$B=I\backslash M$\textit{\ with }$g$ \textit{%
satisfying} (\textit{Res}$CFE$) \textit{on} $B,$ \textit{there is }$\lambda
_{0}\in B,$\textit{\ such that with }$c:=g(\lambda _{0})/\lambda
_{0}^{\kappa },$%
\[
g(\lambda )=c\lambda ^{\kappa }\qquad (\lambda \in B);
\]%
$\ $

\noindent (ii) \textit{With }$g$ \textit{and }$B$\textit{\ as above, }$f,$ $%
\{x_{n}\}_{n\in \mathbb{N}},\{a_{n}\}_{n\in \mathbb{N}}$\ \textit{as in
Kendall's Theorem, i.e. }%
\[
a_{n}f(\lambda x_{n})\rightarrow g(\lambda )\qquad (\lambda \in B),
\]%
$f$\textit{\ is regularly varying with index }$\kappa $:%
\[
f(x)\sim x^{\kappa }\ell (x)
\]%
\textit{with }$\ell $\textit{\ slowly varying.}

\noindent (iii) \textit{In }(ii), \textit{the sequence }$\{a_{n}\}_{n\in
\mathbb{N}}$\textit{\ is regularly varying relative to }$\{x_{n}\}_{n\in
\mathbb{N}}$\textit{\ with index }$-\kappa :$%
\[
a_{n}\sim cx_{n}^{-\kappa }/\ell (x_{n}).
\]

\bigskip

\noindent \textbf{Proof. }(i) We denote by $M(s)$ the exceptional negligible
$\lambda $-set appearing in (\textit{Res}$CFE$)$_{B}.$ W.l.o.g. these are
meagre $\mathcal{F}_{\sigma }$ in the category case.

We begin in (a) below by proving a local version of (i).

(a) Here we take $I:=(a,b)\ $with $0<a<b$ (and $B=I\backslash M,$ with $M$
is negligible). By assumption $g$ is continuous on $B.$ Fix $\varepsilon $
with%
\[
0<\varepsilon <\left( \sqrt{b/a}\right) -1.
\]%
Then $a(1+\varepsilon )<b/(1+\varepsilon ),$ and so the interval%
\[
I_{\varepsilon }:=(a(1+\varepsilon ),b/(1+\varepsilon
))=\tbigcap\nolimits_{s\in J_{\varepsilon }}(I\cap s^{-1}I)\qquad \text{for }%
J_{\varepsilon }:=(1/(1+\varepsilon ),(1+\varepsilon ))
\]%
is non-degenerate. So, for $s\in J_{\varepsilon }$: $I_{\varepsilon
}\subseteq I\cap s^{-1}I,$ in particular, if $\lambda _{0}\in I_{\varepsilon
},$ then $\lambda _{0}\in s^{-1}I$. But (adapting the notation used earlier
for $C(s))$%
\begin{eqnarray*}
B\cap s^{-1}B &=&(I\backslash M)\cap s^{-1}(I\backslash M)=(I\backslash
M)\cap (s^{-1}I\backslash s^{-1}M) \\
&\supseteq &C(s):=I_{\varepsilon }\backslash (M\cup s^{-1}M).
\end{eqnarray*}

By (\textit{Res}$CFE$)$_{B}$, for $\varepsilon $ small enough and all $s\in
J_{\varepsilon }$, for quasi all $\lambda \in C(s)$%
\[
s^{\kappa }=K(s)=g(s\lambda )/g(\lambda ).
\]%
Let $D=\{d_{n}\}_{n\in \mathbb{N}}$ enumerate a countable set dense in $%
J_{\varepsilon }.$ The sets $d_{m}^{-1}M,d_{m}^{-1}M(d_{m})$ being
negligible (meagre $\mathcal{F}_{\sigma }$, in the category case), the set%
\[
H:=\bigcap\nolimits_{m}I_{\varepsilon }\backslash \lbrack d_{m}^{-1}M\cup
d_{m}^{-1}M(d_{m})]\subseteq I,
\]%
is non-negligible and so non-empty (in the category case: a dense $\mathcal{G%
}_{\delta }$ in $I_{\varepsilon }$, so the Baire Category Theorem applies).
Take $\lambda _{0}\in H.$

As above, $\lambda _{0}\in d_{m}^{-1}I,$ since $\lambda _{0}\in
I_{\varepsilon }$, i.e. $\lambda _{0}d_{m}\in I.$ Also $\lambda
_{0}d_{m}\notin M,$ i.e. $\lambda _{0}d_{m}\in I\backslash M.$ So $g$ is
continuous at each $\lambda _{0}d_{m}.$ Furthermore, $\lambda _{0}d_{m}\in
I\backslash M(d_{m})$, so
\[
d_{m}^{\kappa }=g(\lambda _{0}d_{m})/g(\lambda _{0}):\qquad g(\lambda
_{0}d_{m})=d_{m}^{\kappa }g(\lambda _{0}).
\]%
But $\{d_{m}\}_{m}$ is dense in the interval $J_{\varepsilon },$ and also,
as we have seen, $\lambda _{0}d_{m}\in I\backslash M.$ So, by continuity of $%
g$ on $B,$ passage to the limit gives, for all $\lambda _{0}t$ in $\lambda
_{0}J_{\varepsilon }\cap B,$ i.e. for quasi all $t$ in $\lambda
_{0}^{-1}J_{\varepsilon },$ that
\[
g(\lambda _{0}t)=t^{\kappa }g(\lambda _{0}).
\]%
Writing $\lambda =\lambda _{0}t,$ for quasi all $\lambda \in \lambda
_{0}J_{\varepsilon }\cap B,$%
\[
g(\lambda )=\lambda ^{\kappa }g(\lambda _{0})\lambda _{0}^{-\kappa }.
\]

(b) The argument in (a) above may be repeated, mutatis mutandis, in any
subinterval of $I,$ and this will allow us to prove (i). We put $%
h(x):=g(x)x^{-\kappa }$ and%
\[
J:=\{x\in B:(\exists k_{x})(\exists \delta )\text{ }h|(xJ_{\delta
})=k_{x}\}.
\]%
Then $J$ is open in $B$, and by the earlier argument everywhere dense in $I.$
Consider any maximal interval $J^{\prime }:=(a^{\prime },b^{\prime })$ with $%
J^{\prime }\cap B$ contained in $J$ and let $k=h|(J^{\prime }\cap B)$.
Suppose that $b^{\prime }$ is interior to $I.$ For $s$ with $s<1$ the
interval $s^{-1}J^{\prime }$ contains $b^{\prime }(1,s^{-1})$ and so meets
all the maximal intervals $(c,d)$ of $J$ sufficiently close on the right of $%
b^{\prime }.$ Fix any such interval $(c,d).$ So $s^{-1}B\cap B\supseteq $ $%
s^{-1}J\cap J\supseteq s^{-1}(J^{\prime }\cap B)\cap (c,d).$ Select $%
s\lambda \in J^{\prime }\cap B$ with $\lambda =s^{-1}(s\lambda )\in
(c,d)\cap B\backslash M(s)$. Let $k^{\prime }=h|(c,d)\cap B.$ Then as $%
\lambda \in (B\cap s^{-1}B)\backslash M(s)$%
\[
g(s\lambda )=s^{\kappa }g(\lambda ):\qquad k=h(s\lambda )=h(\lambda
)=k^{\prime }.
\]%
Thus, on any maximal interval sufficiently close to $b^{\prime },$ $h=k;$
this contradicts the maximality of $(a^{\prime },b^{\prime })$ unless $%
b^{\prime }$ is not an interior point of $I.$ So $b^{\prime }=b$ the right
end-point of $I$. Likewise, $a^{\prime }=a$ the left end-point of $I.$ So $%
h(x)$ is constant on $B.$ $\square ($i$).$

(ii) In the proof of Theorem 1 in \S 2 we showed that for some $\kappa $ the
function $g$ satisfies (\textit{Res}$CFE$) on $B$ and that $f$ is regularly
varying with index that $\kappa .$ $\square ($ii$).$

(iii) By (ii) write $f(x)\sim x^{\kappa }\ell (x),$ for some $\ell $ $\in SV$%
, and by (i) $g(\lambda )=c\lambda ^{k}$ for $\lambda \in B.$ So, for any
fixed $\lambda \in B,$
\[
a_{n}\sim g(\lambda )/f(\lambda x_{n})=c\lambda ^{\kappa }/[\lambda ^{\kappa
}x_{n}^{\kappa }\ell (\lambda x_{n})]\sim cx_{n}^{-\kappa }/\ell
(x_{n}).\qquad \square
\]

\noindent \textbf{Corollary. }\textit{For }$B\ $\textit{Baire and }$g$
\textit{satisfying} (\textit{Res}$CFE$) \textit{on} $B,$ \textit{there is a
discrete family of intervals }$I$\textit{\ each with a corresponding point }$%
\lambda _{I}\in I,$\textit{\ such that with }$c_{I}:=g(\lambda _{I})/\lambda
_{I}^{\kappa },$%
\[
g(\lambda )=c_{I}\lambda ^{\kappa }\qquad (\lambda \in I).
\]

\bigskip

\noindent \textbf{Proof. }By Theorem 6I, the quasi-interior $B^{q}$ of $B$
may be represented as a union of open intervals on each of which $h$ is
constant. Consider the family of maximal open intervals in $B^{q}$ on which $%
h$ is constant. Then, as in the proof of (i), both end points of such a
maximal interval are not limits of other maximal intervals. $\square $

\bigskip

So far we have considered the general Baire and the interval-minus-null
cases. We turn now to the general measure variant: this addresses the
measure case, and identifies a scenario of $h$-constancy on certain
`rational skeletons', dependent on the points of $B.$ See [BinO1, Th. 4.2]
for the constancy of rationally invariant functions (i.e. when $h(qx)=h(x)$
for all $q\in \mathbb{Q}).$ Here differences arise between Theorem 6M and 6B
because Theorem 6M requires \textit{quantitative} measure theory (rather
than \textit{qualitative} measure theory, which is closely aligned with the
Baire case). The breakdown of the usual category-duality occurs here since a
measurable set need not be `locally co-null' at any point (i.e. never meets
almost all of some interval); the analogous qualitative argument delivers
less information.

In the Corollary above, one has multiple constancy. Theorem 6M\ below in the
case when $B$ is non-null and nowhere dense opens a similar possibility.

\bigskip

\noindent \textbf{Theorem 6M.}\textit{\ Take }$B$\textit{\ non-null closed
and }$g$ \textit{satisfying} (\textit{Res}$CFE$) \textit{on} $B.$

\noindent (a) \textit{If }$B$\textit{\ is nowhere dense, then for almost all
}$b\in B$\textit{\ there exists in }$\mathbb{Q}$ \textit{a sequence\ }$%
q_{n}=q_{n}(b)\rightarrow 1$\textit{\ with }$q_{n}b\in B$\textit{\ and}%
\[
g(q_{n}b)=q_{n}^{\kappa }g(b):\qquad g(bq_{n}(b))/(bq_{n}(b))^{\kappa
}=g(b)/b^{\kappa },
\]%
\textit{i.e. }$h(x):=g(x)/x^{\kappa }$ \textit{remains constant on a
rational sequence of dilations }$q_{n}(b)$\textit{. The sequence can be
selected so that }$q_{2n-1}(b)\downarrow 1$ \textit{and} $%
q_{2n}(b)=1/q_{2n-1}(b)\uparrow 1.$

\noindent (b) \textit{If }$B$\textit{\ is somewhere dense, then }$%
B=B_{1}\cup B_{0}$ \textit{with }$B_{1}$ \textit{open} \textit{and} $B_{0}$
\textit{nowhere dense; then Theorem 6B(i)-(iii) applies\ mutatis mutandis to
}$B_{1},$\textit{\ and (a) above applies to }$B_{0}$\textit{\ if }$B_{0}$
\textit{is non-null.}

\bigskip

\noindent \textbf{Proof.} We follow the notation of Th. 6I, in particular,
we write $J_{\Delta }:=[(1+\Delta )^{-1},(1+\Delta )].$

(a) W.l.o.g. $B$ is (closed and) of finite measure and $g$ is continuous on $%
B.$ Now notice that the set%
\[
S(B):=\{\lambda \in B:(\exists \{s_{n}\}\uparrow 1)[g(s_{n}\lambda
)=s_{n}^{\kappa }g(\lambda )\text{\&(}\lambda \in s_{n}^{-1}B)]\}
\]%
is analytic, hence measurable. We will show that $S(B)$ is non-empty, and
hence is almost all of $B$: indeed, suppose otherwise, then $B\backslash
S(B) $ is non-null. Then passing to a non-null closed subset, $F\ $say, it
follows that $\emptyset \neq S(F)\subseteq F\cap S(B),$ contradicting that $%
F $ is disjoint from $S(B).$

Notice that density-open subsets meet in a density-open subset (empty or
otherwise).

Fix $p$ with $1/2<p<1$ (e.g. $p=3/4).$ We now define an operator $\Delta $
on on-empty density-open sets $A\subseteq B.$ Choose a density point $b\in
A. $ There is $\Delta =\Delta (A)>0$ such that%
\[
|A\cap bJ_{\Delta }|>p|bJ_{\Delta }|.
\]%
So taking $q=1-p$%
\[
|bJ_{\Delta }\backslash A|<q|bJ_{\Delta }|.
\]%
For $s>1$ with $s\in J_{\Delta }$ one has $s-1<\Delta $ $,$ and $%
0<1-s^{-1}<1-1/(1+\Delta )=\Delta /(1+\Delta ).$ Take $L(\Delta ):=\max
\{(1+\Delta ),\Delta /(1+\Delta )\}=1+\Delta ,$ so that $L(\Delta
)\rightarrow 1$ as $\Delta \rightarrow 0.$ Now for $s\in J_{\Delta }$
\[
|bJ_{\Delta }\backslash bJ_{\Delta }s^{-1}|\leq |1-s|bL(\Delta ).
\]%
So for $s\in J_{\Delta }$,
\[
|As^{-1}\cap A\cap (bJ_{\Delta }\backslash A)|>(p-q)|bJ_{\Delta
}|-|1-s|bL(\Delta ).
\]

So in the nhd $bJ_{\Delta }$ of $b,$ for any $s$ close enough to $1,$ and on
either side of $1,$ the set $C(A):=(As^{-1}\cap A\cap bJ_{\Delta })$ is
density-open, qua intersection of density-open sets. As it has non-null
measure it is non-empty, so contains a density point, $c(A)$ say.

We work inductively. Base step: taking $C_{0}$ to be the density-open
interior of $B$ (with $b\in C_{0}),$ put $\Delta _{1}:=\Delta (C_{0}),$ with
$\Delta \ $the operator above. Take $s_{1}\in J_{\Delta (1)}$ small enough,
as above, so that $C_{1}:=C(C_{0})$ is density-open and contains $b.$

Continue selecting $s_{1},s_{2},...\rightarrow 1,$ with $s_{2i}$ increasing
and $s_{2i+1}$ decreasing and rational, and non-empty density open sets $%
C_{i}\subseteq B$ with distinguished member $b_{i}=c(C_{i}),$ so that with $%
\Delta (i):=\Delta (C_{i}):$

i) $C_{i+1}:=C(C_{i})=(C_{i}s_{i}^{-1}\cap C_{i}\cap b_{i}J_{\Delta (i)})$
is non-empty and density-open so contains a density point $b_{i+1}$

ii) $\lambda (s_{i})\in C_{i}\backslash \tbigcup\nolimits_{j\leq i}M(s_{j})$.

Then $\lambda (s_{i})\in (C_{i}s_{i}^{-1}\cap C_{i}\cap b_{i}J_{\Delta
(i)})\subseteq (C_{j}s_{j}^{-1}\cap C_{j}\cap bJ_{\Delta (i)})\subseteq B$
and so $\lambda (s_{i})s_{j}\in C_{j}\subseteq B.$ Thus $\lambda (s_{i})\in
B\cap Bs_{j}^{-1}.$ To apply (\textit{Res}$CFE)_{B},$ note that also $%
\lambda (s_{i})\notin M(s_{j});$ so
\[
g(s_{j}\lambda (s_{i}))=s_{j}^{\kappa }g(\lambda (s_{i})).
\]%
But $\lambda (s_{i})\in B,$ so $\lambda _{0}=\lim_{n}\lambda (s_{n}),$ and $%
\lambda _{0}s_{j}=\lim_{n}\lambda (s_{n})s_{j}\in B$, as $B$ is closed. By
continuity of $g$%
\[
g(s_{j}\lambda _{0})=s_{j}^{\kappa }g(\lambda _{0}).
\]

(b) If $B$ is (closed and) somewhere dense, take $B_{1}$ the union of
maximal open intervals contained in $B.$ Then $B_{0}:=B\backslash B_{1}$ is
closed and nowhere dense (otherwise, it would contain an open interval
disjoint from $B_{1}).$ The remaining assertions are clear. $\square $

\bigskip

The results above extend to general regular variation. The proof is much as
above, via the Popa circle groups, with $s^{\kappa }=g(\lambda s)/g(\lambda
) $ and $(B\cap s^{-1}B)$ above replaced by $K(s)=g(\lambda \circ _{\eta
}s)/g(\lambda )$ and $(B\cap s_{\eta }^{-1}\circ _{\eta }B),$ respectively;
we omit the details.

\bigskip

\noindent \textbf{Remark. }The key idea of Theorem 6 is the embedding of a
specific countable set, one that is dense in itself, into an open set
`punctured' by the removal of a small or negligible part. Embeddings of
countable sets by translation go back to Marczewski [Mar]; see [NatO] for
recent developments.

\bigskip

\textbf{8. Character degradation from ess-lim.}

8.1 In what follows, we will need to distinguish between (general) sets of
reals, and sets which can be defined by suitable coding. For background
here, see e.g. the monograph Kechris [Kec, Ch.V] on the analytical hierarchy
(note [Kec, V.40B] on classical v. effective descriptive set theory), and
our recent survey [BinO9]. For a deeper analysis of coding see [Sol, II.1.1,
25-33]; a minimal amount is in [FenN, \S\ 2, p. 93]. We defer further
discussion of these matters (including the ambiguous analytical class $%
\mathbf{\Delta }_{2}^{1})$ to the Appendix below, and to the proof of
Theorem 7 below.

To say that $L=$ess-$\lim_{x\rightarrow \infty }f(x)$ requires the assertion
that \textit{there exists} a (meagre, exceptional)\textit{\ set} off which
\textit{for all }$x$\textit{\ real }(large enough) $f(x)$ is as close to $L$
as desired. In brief this has an $\exists \forall $ quantifier block in
regard to the `analytical objects': \textit{sets} and \textit{reals}.

It is also true that its negation $L\neq $ess-$\lim_{x\rightarrow \infty
}f(x)$, the assertion that \textit{there exists} a (non-meagre) \textit{set}
on which \textit{for all }$x$ \textit{real }(large enough) $f(x)$ avoids
being sufficiently close to $L,$ also has an $\exists \forall $ structure in
regard to analytical objects.

These two observations have a rigorous formulation below, which adds to
earlier considerations of the character of limits in Karamata and Beurling
RV noted already in [BinO2,5]. Though our proof is largely self-contained,
we refer for background and for the notation of the analytical hierarchy
needed here to [Kec, Ch.V], [BinO9], and to the Appendix below.

\bigskip

\noindent \textbf{Theorem 7 (Character degradation). }\textit{For }$k$%
\textit{\ Borel, the predicate }%
\[
K(s):=\mathrm{ess}\text{-}\lim_{x\rightarrow \infty }k(s,x)
\]%
\textit{is of ambiguous analytical class} $\mathbf{\Delta }_{2}^{1}.$

\bigskip

\noindent \textbf{Proof. }For simplicity, we consider the equation $L=$ess-$%
\lim_{x\rightarrow \infty }f(x)$ with $f$ Baire, say. This is equivalent to
the predicate
\[
(\exists a\in \mathbb{R})(\forall x\in \mathbb{R})(\forall m,n\in \mathbb{N}%
)(\exists p\in \mathbb{N})\Phi (f,a,x,m,n,p),
\]%
where the matrix $\Phi $ is%
\[
\lbrack G(a(n))\text{ is everywhere dense}]\text{\&}[x\in G(a(n))\text{\&}%
x>p\Rightarrow |f(x)-L|<1/m].
\]

Here $\Phi $ `says' that, on the dense $\mathcal{G}_{\delta }$ set $%
\tbigcap\nolimits_{n}G(a(n))$ and to the right of $p,$ the values $f(x)$ are
to within $1/m$ of $L;$ here $a(n):=a\cap \{1\cdot 2^{n},3\cdot 2^{n},5\cdot
2^{n},$ $...\}.$ Apart from the \textit{arithmetic} (= natural number)
quantifiers acting on $\Phi $ there are two quantifiers, the first
existential, the second universal, ranging over the analytical objects of
type 1, the real numbers $a$ and $x$; in view of the opening \textit{%
analytical} quantifier block $\exists \forall $ of $2$ quantifiers over type
1 objects, the statement is said to be $\Sigma _{2}^{1}(f)$ -- the
parentheses acknowledge use of $f$ as an input. The (light-faced, here)
sigma symbol identifies the first quantifier as existential. Under our
simplifying assumption of \S\ 2 (Remark 2) that a Baire $f$ is replaced by a
`nice' Borel version with all its rational level sets being $\mathcal{G}%
_{\delta }$, we can think of $\Phi $ as written not with the use of $f$ but
instead in terms of two codes (= real numbers) $b$ and $c,$ where $%
b(m)=b_{m} $ and $c(m)=c_{m}$ as above, for each $m.$ For details see the
Appendix. To indicate a hidden mention of the need for some real parameters,
we use a bold symbol,\ and say more simply that the statement is $\mathbf{%
\Sigma }_{2}^{1}.$

As noted earlier, the negation can also be expressed as a $\mathbf{\Sigma }%
_{2}^{1}$ statement: the inequality $L\neq $ess-$\lim_{x\rightarrow \infty
}f(x)$ is equivalent to a statement of the form
\[
(\exists a\in \mathbb{R})(\forall x\in \mathbb{R})(\exists m\in \mathbb{N}%
)(\forall p\in \mathbb{N})\Psi (a,x,m,p),
\]%
where the matrix $\Psi $ is the statement%
\[
\lbrack G(a(m))\text{ unbounded \&}[[x\in G(a(m))\&(x>p)]\Rightarrow
|f(x)-L|\geq 1/m]].
\]%
We summarize the above by saying that ess-lim is of \textit{ambiguous
analytical class} $\mathbf{\Delta }_{2}^{1}.$

The analysis above extends with little extra complexity to cover the case $%
L(s)=$ess-$\lim_{x\rightarrow \infty }k(s,x).$ $\square $

\bigskip

\noindent \textbf{Remark. }Character degradation under ess-lim here amounts
(for $f$ Borel) to the Borel set
\[
H_{m}(f):=\{x:|k(x)-L|<1/m\}
\]%
rising to the second level of the analytical hierarchy and becoming a more
complex set: a $\mathbf{\Delta }_{2}^{1}$ set. (With regard to the set $%
H_{m}(f)$ see the Appendix.)

\bigskip

8.2. \textit{Provably }$\Delta _{2}^{1}$\textit{\ sets. }A set $A$ is said
to be \textit{provably} $\Delta _{2}^{1}$ if there are two $\Sigma _{2}^{1}$
predicates $\Phi (x,y)$ and $\Psi (x,y)$ and a real number $b$ such that $%
a\in A$ iff $\Phi (a,b),$ and likewise $a\notin A$ iff $\Psi (a,b),$ with
both these equivalences provable in the axiom system $ZF+DC,$ where $DC$
stands for the Axiom of Dependent Choices. Thus $A$ is in $\mathbf{\Delta }%
_{2}^{1}$ (bold-faced, because of the parameter $b$).

Fenstad and Normann [FenN] noticed that a key step in [Sol] (in which
Solovay constructs a model of set theory wherein $DC$ holds and all subsets
of $\mathbb{R}$ are Lebesgue measurable) may be re-read to show that all
\textit{provably} $\Delta _{2}^{1}$ subsets are measurable: see Remarks 1
and 3 in [FenN, p. 95]. Ultimately, the argument relies on the notion of
forcing provided by the partially ordered set comprising the Borel sets
lying in a fixed countable model of set theory -- see [BinO9, \S\ 6.1]. By
using the category variant of this partial order, much the same argument
gives that all \textit{provably} $\Delta _{2}^{1}$ subsets have the Baire
property: see [Kan, \S\ 14.4, p.180].

In conclusion: we should not be surprised that `nice versions' of Baire
functions yield corresponding essential-limit functions that are Baire.

\bigskip

\noindent \textbf{9. Complements.}

\noindent 9.1 \textit{Axiomatics: set-theoretic foundations for RV. }We have
stressed in the Introduction the role of the Axiom of Dependent Choice(s),
DC. Its great strength, as Solovay [Sol, p. 25] points out, is that it is
sufficient for the establishment of Lebesgue measure, i.e. including its
translation invariance and countable additivity ("...positive results ... of
measure theory..."), and may be assumed consistently with such additional
axioms as LM (\textit{all} subsets are Lebesgue measurable) and PB (\textit{%
all} subsets have Baire property, BP). To generate non-measurable sets one
needs the Axiom of Choice AC. While the Zermelo-Fraenkel(-Skolem) axiom
system ZF is common ground in mathematics, AC is not, and alternatives to it
are widely used, including the two we have just mentioned. For a thorough
discussion of alternatives we refer to [BinO9], especially \S\ 10 therein.

In the standard Karamata setting of RV, continuous limits of functions may
be replaced by sequential limits (as in \S\ 2 above), so that starting with
continuous functions one remains within the class of Borel functions. In
replacing limits by limsups character degradation occurs leading functions
up the analytical (projective) hierarchy. In this connection we have
previously argued [BinO2, \S\ 5] that $\mathbf{\Delta }_{2}^{1}$ is a most
attractive class of sets within which to carry out the analyses of RV. When
drawing in the Beurling operation $\circ _{\varphi },$ this argument needed
amplification -- see [BinO6, \S\ 11]. Here, when extending the argument to
`essential limits' we again point to the further attractions of the\textit{\
provably} $\mathbf{\Delta }_{2}^{1}$ class of \S\ 6.2. Working with `nice
versions' of functions -- removing a pathological set covered by a $\mathcal{%
G}_{\delta }$ set, as in the Appendix -- via the Baire-Kuratowski or Luzin
Continuity Theorems [Oxt, Th. 8.1, 8.2], we remain in the realm of
Baire/measurable sets, as though under the sway of LM or PB. This is because
the Baire Category Theorem, BC, suffices here. Indeed, BC is equivalent to
DC; see [BinO9] and the literature cited there.

\noindent 9.2 \textit{Smallwood's theorem. }Essential (or approximate)
limits go back to work of Denjoy in 1916 on approximate continuity, and
Khintchine in 1924 and 1927. An early textbook treatment is in Saks [Sak,
IX.10]. Smallwood's Theorem [Sma] reconciles the Denjoy and Khintchine
approaches: for $E\subseteq \mathbb{R}$ Lebesgue-measurable, $f:\mathbb{%
R\mapsto R}$ measurable, $x_{0}\in \mathbb{R},$ $f$ has approximate limit $L$
at $x_{0}$ (i.e. $\exists $ess- lim$f(x_{0})=L$) if and only if there exists
a measurable set $F\subseteq E$ with density 1 at $x_{0}$%
\[
\lim\nolimits_{x\in F,x\rightarrow x_{0}}f(x)=L
\]

\noindent Such matters are important in probability theory; see e.g. the
survey of Geman and Horowitz [GemH, Appendix: Metric density and approximate
limits, 22-24], and the lecture notes of Adler [Adl, IV.4.6].

For essential (or approximate) semi-continuity, see Zink [Zin].

\noindent 9.3 \textit{Croftian theory and admissible sequences. }Croft's
theorem says that for a continuous function $f$, the existence of all the
sequential limits (as $n\rightarrow \infty $) of $f(nh)$ for all $h>0$
implies that of the continuous limit of $f(x)$. Kingman [Kin1,2] re-writes
this additively, so working with $f(\log n+x)$, and generalizes the Croft
setting to ask for conditions on $c_{n}$ for a similar result to apply for
sequential limits if $f(c_{n}+x)$. Roughly speaking, the condition needed
for a croftian theorem to hold here is the Kingman condition%
\[
c_{n+1}-c_{n}\rightarrow 0
\]%
(compare our \textit{admissible sequences}, in additive and multiplicative
notations). As in the work of Kingman [Kin1,2] and Kendall [Ken], the key
role of the Baire category theorem is clear in the following further
generalization of Vinokurov [Vin]: for $f$ as above and $c_{n}$ satisfying
Kingman's condition, the condition on the set $E$ needed for the implication
from%
\[
f(c_{n}+x)\rightarrow L\qquad (x\in E)
\]%
to%
\[
f(x)\rightarrow L
\]%
is that $E$ be non-meagre. For further results of this type, see Feh\'{e}r
at al. [FehLT]; cf. [Sen].

\noindent 9.4 \textit{Regularly varying measures. }For random vectors $X$ in
$\mathbb{R}^{d}$, a theory of regularly varying measures can be based on the
definition (suggested by Kendall's Theorem)
\[
n\mathbb{P}(X/a_{n}\in \text{ }.)\rightarrow \mu (.)\qquad (n\rightarrow
\infty )
\]%
for $a_{n}\nearrow \infty $, vague convergence, and suitably restricted $\mu
$. Then regular variation is present, as for some $\alpha >0$%
\[
\mu (tA)\sim t^{\alpha }\mu (A)
\]%
as in Theorem 6 (and then $a_{n}$ is regularly varying). See e.g. Hult and
Linskog [HulL], Hult et al. [HulLMS]. This approach is Kendall-like, as it
is entirely sequential. It can be extended to infinite-dimensional settings,
and is widely used nowadays in probability (theory and applications).

\noindent 9.5 \textit{Thinning:\ Steinhaus-Weil aspects. }The two main
ingredients in verifying that the Kendall criterion yields regular variation
are: the croftian property of the set $C$ in Prop. 2, and the Steinhaus-Weil
property of the test set $B$ -- that $BB^{-1}$ contains an interval $J$
around $1.$ Recall that the latter guarantees that $C(s)=B\cap s^{-1}B$ is
non-empty for $s\in J$, and so for $s\in J\ $and $\lambda \in C(s)$
\[
k_{s}(\lambda x_{n})\rightarrow g(s\lambda )/g(\lambda ),\text{ as }%
n\rightarrow \infty .
\]%
The Steinhaus-Weil property can hold for nowhere dense sets (cf. a
multiplicative analogue of the classical Cantor excluded middle thirds);
indeed there is a rich family of such sets -- see the SW\ property used in
[BinO8].

However, Prop. 2 relies both on the Baire function $g$ having a point of
continuity in $C,$ and on $C$ having the property that, for an additively
admissible sequence $c_{n},$ the tail union of its translates $%
\tbigcup\nolimits_{n\geq m}(c_{n}+C)$ contains quasi all of an infinite
half-ray. Just as in Vinokurov's result (\S\ 9.3), $C\ $here cannot be
negligible.

That said, we note that, taken together, Theorem 1 and 1M already imply (on
taking $a_{n}:=1/f(x_{n}))$ the Characterization Theorem K of \S\ 2 with the
global hypothesis of convergence $f(tx)/f(t)\rightarrow g(t)$ `for all $t$'
much weakened to `for all $t$ on a \textit{non-negligible} \textit{set}'. On
the other hand, it is known [BinO8] that a further thinning is possible: to
sets having the SW property locally. One is thus led to ask, given the
convergence \textrm{ess-lim}$_{x\rightarrow \infty }k_{s}(x)$ on an $s$%
-interval, whether the Theorem 1 can be further thinned, despite the
possible negligibility of $C(s)$ itself.

\noindent 9.6 \textit{Functional equations. }For an account of the
literature of GS and related equations see [AczD, Ch. 19] and the more
recent [Brz], cf. the summary in [Ost2, \S\ 1]. In our context it is natural
to restrict solutions of ($GS$) and of the related Chudziak-Jab\l o\'{n}ska
equation
\begin{equation}
H(x\circ _{\eta }y)=H(x)H(y),  \tag{$CJ$}
\end{equation}%
with $\eta $ continuous, to be non-negative and locally bounded. It then
emerges from [BrzM] (cf. [Ost4, \S\ 9.5] for a more direct approach),
[Jab1-4] and especially [Jab3], that, provided the function $H$ is \textit{%
non-trivial} (i.e. its range is not a subset of $\{0,1\}$), then local
boundedness of the solution $H$ implies continuity. (Note the trivial
counter-example: the Dirichlet function $H=\mathbf{1}_{\mathbb{Q}}$ for $%
\eta (t)=1+t.$) This observation includes the case of solutions of ($GS$)
which take the form $\eta (t)=\eta _{\rho }(t)$ for some $\rho \geq 0$ and $%
t>-1/\rho .$ The case $\rho =\infty ,$ corresponding to $x\circ _{\eta
}y=xy, $ is just another instance of ($CFE$).

By Theorem 6, taking $g$ non-trivial in Theorems 4 and 4M ensures the
corresponding $K\ $is likewise non-trivial and so continuous.

Matters are the same in the more general ($BG$) equation. The case for $%
\sigma (t)=1+st$, with $s\geq 0$ is typified by $s=1$ via scaling (save for
the case $s=0,$ reduced via logarithms to $s=1);$ then $u\circ _{\sigma
}v=u+v+uv=(u+1)(v+1)-1$ and here the ($BG$) equation reduces to%
\[
K(x\circ _{\eta }y)+1=(K(x)+1)(K(y)+1),
\]%
so that $H(x):=K(x)+1$ is locally bounded and so continuous provided $K$
(being non-negative) is non-zero -- by Jab\l o\'{n}ska's theorems in [Jab3].
But here, again Theorem 6, since $g$ is assumed non-zero in Theorems 5 and
5M, ensures the corresponding $K\ $is again continuous.

The continuous solutions of ($BG$) are given in the table below. (The FE
literature also includes studies of the case where on the right $\circ
_{\sigma }$ is replaced by a semigroup operation $\circ $ as in [Chu1,2].)

In the table, the four corner-formulas correspond to classical variants of
the Cauchy functional equation ($CFE$)$.$ For completeness we include the
proof; this proceeds by a straightforward reduction to a classical variant
of ($CFE$) by an appropriate shift and rescaling, similar to the reduction
from $K$ to $H$ above. The notation $\circ _{r}$ etc. below refers to the
Popa operation $\circ _{\eta }$ with parameter $r,$ i.e. the case $\eta
=\eta _{r}.$

\bigskip

\noindent \textbf{Proposition 3 }([Ost3, Prop. A; cf. [Chu1])\textbf{. }%
\textit{For }$\circ _{\eta }=\circ _{r},\circ _{\sigma }=\circ _{s},$\textit{%
\ and }$K$ \textit{Baire/measurable satisfying }($BG$)$,$ \textit{there is }$%
\kappa \in \mathbb{R}$ \textit{so that }$K(t)$ \textit{is given by:}%
\renewcommand{\arraystretch}{1.25}%
\[
\begin{tabular}{|l|l|l|l|}
\hline
Popa parameter & $s=0$ & $s\in (0,\infty )$ & $s=\infty $ \\ \hline
$r=0$ & $\kappa t$ & $(e^{\kappa t}-1)/s$ & $e^{\kappa t}$ \\ \hline
$r\in (0,\infty )$ & $\kappa \log (1+rt)$ & $[(1+rt)^{\kappa }-1]/s$ & $%
(1+rt)^{\kappa }$ \\ \hline
$r=\infty $ & $\kappa \log t$ & $(t^{\kappa }-1)/s$ & $t^{\kappa }$ \\ \hline
\end{tabular}%
\]%
\renewcommand{\arraystretch}{1}

\noindent \textbf{Proof. }Each case reduces to ($CFE$) on $\mathbb{R}_{+}$,
or a classical variant by an appropriate shift and rescaling. For instance,
given $K,$ for $r,s>0$ set%
\[
F(t):=1+sK((t-1)/r):\qquad f(\tau )=(K(1+r\tau )-1)/s.
\]%
Then with $u=1+rx,v=1+ry,$ as $(uv-1)/r=x\circ _{r}y,$%
\[
F(uv)=1+sK(x\circ _{r}y)=1+sK(x)+sK(y)+s^{2}K(x)K(y)=F(u)F(v),
\]%
for $u,v\geq 0.$ So, as $F$ is Baire/measurable (see again [Kuc3, \S\ 13]), $%
F(t)=t^{\gamma }$ and so $K(t)=[(1+rt)^{\gamma }-1]/s.$ The remaining cases
are similar. $\square $

\bigskip

In the language of isomorphisms $\eta _{\rho },$ $\exp ,$ $\log $, we can
rephrase the above more succinctly as follows:

\renewcommand{\arraystretch}{1.25}%
\[
\begin{tabular}{|l|l|l|l|}
\hline
Popa parameter & $\sigma =0$ & $\sigma \in (0,\infty )$ & $\sigma =\infty $
\\ \hline
$\rho =0$ & $\kappa t$ & $\eta _{\sigma }^{-1}(e^{\kappa t})$ & $e^{\kappa t}
$ \\ \hline
$\rho \in (0,\infty )$ & $\log \eta _{\rho }(t)^{\kappa }$ & $\eta _{\sigma
}^{-1}(\eta _{\rho }(t)^{\kappa })$ & $\eta _{\rho }(t)^{\kappa }$ \\ \hline
$\rho =\infty $ & $\log t^{\kappa }$ & $\eta _{\sigma }^{-1}(t^{\kappa })$ &
$t^{\kappa }$ \\ \hline
\end{tabular}%
\]%
\renewcommand{\arraystretch}{1}

\noindent 9.7 \textit{Open Question.}\textbf{\ }In passing, motivated by the
context of Prop. 3, we leave open the question whether Prop. 1 (based on an
affine action between two sets $A,$ $B$) has an analogue for more general
continuous $h(z,s)$ in the spirit of the `Miller homotopies' [Mil], cf.
[BinO4]. (We have in mind something along the lines: for convergent
sequences $z_{n}\rightarrow z_{0},$ for almost all $b$ near $b_{0}$ there
are infinitely many $m$ with $H(b,z_{m})\in A$ -- here with $H$ the
appropriate local inverse at $b_{0}=h(z_{0},a_{0})$.) See [Gro] as to a
possible approach for replacing differentiability by `Radon-Nikodym
differentiability' as in the `Functionwise Steinhaus-Weil Theorem', wherein $%
f(A\times B)$ has an interior point (originating with Marcin Kuczma [Kuc1]).

\bigskip

\noindent \textbf{Appendix: Relevant aspects of coding.}

Theorem 7 above relied on the ability to refer to various subsets of the
real line, especially open sets, in terms of `codes'. Our canonical sources
there were [Kec, Ch.V] on the analytical hierarchy (and the note [Kec,
V.40B] on classical versus effective descriptive set theory), and our recent
survey [BinO9], and for coding the wide-ranging use in [Sol, II.1.1, 25-33]
and the much more minimal amount in [FenN, \S\ 2, p. 93]. Here we give some
examples to help clarify the full effect on the \textit{analytical
quantifier blocks} in Theorem 7. We begin with some notation.

Let $\{I_{n}\}_{n\in \mathbb{N}}$ enumerate (constructively) all the
rational-ended intervals, with $I_{n}=(l_{n},r_{n})$. Write $\mathbb{M}$ for
the odd natural numbers; for $a\subseteq \mathbb{N}$ we may extract an $n^{%
\text{th}}$ canonical subset of $a$ and also an open set naturally `coded'
by $a$ by setting:%
\[
a(n)=a\cap \{2^{n}m:m\in \mathbb{M}\},\qquad G(a):=\bigcup\nolimits_{n\in
a}I_{n}.
\]%
We identify $a\subseteq \mathbb{N}$ with the real number in $\{0,1\}^{%
\mathbb{N}}$ whose binary expansion is the indicator function of $a$. Thus $%
\{a:m\in a\}$ is open (being the set of reals with $m$-th binary digit =1).

The following are examples of Borel sets, using semi-formal predicates:%
\begin{eqnarray*}
\{a &:&G(a)\text{ is unbounded}\}\Leftrightarrow \\
\{a &:&(\forall k\in \mathbb{N})(\exists q\in \mathbb{Q})(\exists m\in
\mathbb{N})[(q>k)\text{ \& }m\in a\text{ \& }q\in I_{m}]\}, \\
\{a &:&G(a)\text{ is everywhere dense}\}\Leftrightarrow \\
\{a &:&(\forall n\in \mathbb{N})(\exists m\in \mathbb{N})(\exists q\in
\mathbb{Q})[m\in a\text{ \& }q\in I_{n}\text{ \& }q\in I_{m}]\}.
\end{eqnarray*}%
The defining statements here are said to be \textit{arithmetic} since the
quantifiers are `arithmetic': they all range over the (countable) set of
natural or rational numbers, and the `matrix' (the expression in square
brackets -- not containing quantifiers) is built from elementary relations
like $l_{n}<q<r_{n}$ and $m\in a$, and these may be viewed as (codes for the
very simple, basic) open sets containing the real numbers $a$. The two
examples above are Borel, as they may be constructed using countable unions
and intersections (corresponding to the arithmetic quantifiers) from basic
open sets; for example, the second set is $\mathcal{G}_{\delta }$, since it
has the form
\[
\bigcap\nolimits_{n\in \mathbb{N}}\bigcup\nolimits_{m\in \mathbb{N},q\in
\mathbb{Q}}\{a:m\in a\text{ \& }q\in I_{m}\text{ \& }q\in I_{n}\}.
\]%
For $f$ a Baire function and $L$ fixed, the set%
\[
H_{m}(f):=\{x:|f(x)-L|<1/m\}
\]%
is Baire. If $f$ happened to be continuous, this set would be open, and so
coded as $G(a_{m})$ with $a_{m}:=\{n:I_{n}\subseteq H_{m}(f)\},$ i.e. by a
set of natural numbers, or, equivalently, by a single real number. For a
general Baire $f,$ since we are prepared to neglect meagre sets, as
suggested in \S\ 2 Remark 2, we can make a \textit{simplifying assumption}:
regard $H_{m}$ as coded by some open set, $G(b_{m})$ say, less a union of
closed nowhere dense sets -- in essence use a nice version of $f$; passing
to the sequence of the complements of the closed nowhere dense sets, we view
the removal of their union as an intersection of open sets, the $n^{\text{th}%
}$ one coded by the subset $c_{m}(n)$ of $c_{m}$, say. So, for example, with
our simplifying assumption, $H_{m}(f)$ may be regarded as being the $%
\mathcal{G}_{\delta }$ set%
\[
G(b_{m})\cap \bigcap\nolimits_{n\in \mathbb{N}}G(c_{m}(n)).
\]

\noindent \textbf{Acknowledgement. }We thank Tomasz Natkaniec for his very
useful comments.

\bigskip

\noindent \textbf{References}

\noindent \lbrack AczD] J. Acz\'{e}l and J. Dhombres, \textsl{Functional
equations in several variables. With applications to mathematics,
information theory and to the natural and social sciences.} Encycl. Math.
App.\textbf{\ 31}, Cambridge Univ. Press, 1989.\newline
\noindent \lbrack Adl] R. J. Adler, \textsl{An introduction to continuity,
extrema and related topics for general Gaussian processes}. IMS Lect. Notes
\textbf{12}, Inst. Math. Stat., 1990.

\noindent \lbrack Bin1] N. H. Bingham, Tauberian theorems and the central
limit theorem. \textsl{Ann. Prob.} \textbf{9} (1981), 221-231.

\noindent \lbrack Bin2] N. H. Bingham, Scaling and regular variation.
\textsl{Publ. Inst. Math. Beograd} \textbf{97} (111) (2015), 161-174.

\noindent \lbrack Bin3] N. H. Bingham, Riesz means and Beurling moving
averages. \textsl{Risk and Stochastics} (Ragnar Norberg Memorial Volume, ed.
P. M. Barrieu), Imperial College Press, 2019, to appear (arXiv1502.07494v1).

\noindent \lbrack BinGT] N. H. Bingham, C. M. Goldie and J. L. Teugels,
\textsl{Regular variation}. 2nd ed., Cambridge University Press, 1989 (1st
ed. 1987). \newline
\noindent \lbrack BinO1] N. H. Bingham and A. J. Ostaszewski, Kingman,
category and combinatorics. \textsl{Probability and Mathematical Genetics}
(Sir John Kingman Festschrift, ed. N. H. Bingham and C. M. Goldie), 135-168,
London Math. Soc. Lecture Notes in Mathematics \textbf{378}, CUP, 2010.
\newline
\noindent \lbrack BinO2] N. H. Bingham and A. J. Ostaszewski, Regular
variation without limits. \textsl{J. Math. Anal. Appl.} \textbf{370} (2010),
322-338.\newline
\noindent \lbrack BinO3] N. H. Bingham and A. J. Ostaszewski, Dichotomy and
infinite combinatorics: the theorems of Steinhaus and Ostrowski. \textsl{%
Math. Proc. Camb. Phil. Soc.} \textbf{150} (2011), 1-22. \newline
\noindent \lbrack BinO4] N. H. Bingham and A. J. Ostaszewski, Homotopy and
the Kestelman-Borwein-Ditor theorem. \textsl{Canad. Math. Bull.} \textbf{54}
(2011), 12--20.\newline
\noindent \lbrack BinO5] N. H. Bingham and A. J. Ostaszewski, Beurling slow
and regular variation. \textsl{Trans. London Math. Soc., }\textbf{1} (2014)
29-56\newline
\noindent \lbrack BinO6] N. H. Bingham and A. J. Ostaszewski, Beurling
moving averages and approximate homomorphisms. \textsl{Indag. Math. }\textbf{%
27} (2016), 601-633 (fuller version: arXiv1407.4093).\newline
\noindent \lbrack BinO7] N. H. Bingham and A. J. Ostaszewski, Beyond
Lebesgue and Baire IV: Density topologies and a converse Steinhaus-Weil
theorem. \textsl{Topology Appl.}, {\textbf{239} (2018), 274-292 (}%
arXiv1607.00031).\newline
\noindent \lbrack BinO8]{\ N. H. Bingham and A. J. Ostaszewski, Additivity,
subadditivity and linearity: Automatic continuity and quantifier weakening.
\textsl{Indag. Math.} (N.S.) \textbf{29} (2018), 687--713 (arXiv
1405.3948v3).}\newline
\noindent \lbrack BinO9] N. H. Bingham and A. J. Ostaszewski, Set theory and
the analyst. \textsl{European J. Math.}, Online First,
doi.org/10.1007/s40879-018-0278-1 (arXiv:1801.09149v2).\newline
\noindent \lbrack BinO10] N. H. Bingham and A. J. Ostaszewski, General
regular variation, Popa groups, and quantifier weakening. arXiv: 1901.05996.%
\newline
\noindent \lbrack Brz] J. Brzd\k{e}k, The Go\l \k{a}b-Schinzel equation and
its generalizations. \textsl{Aequat. Math.} \textbf{70} (2005), 14-24.%
\newline
\noindent \lbrack BrzM] J. Brzd\k{e}k and A. Mure\'{n}ko, On a conditional Go%
\l \k{a}b-Schinzel equation. \textsl{Arch. Math.} \textbf{84} (2005),
503-511.\newline
\noindent \lbrack Chu1] J. Chudziak, Semigroup-valued solutions of the Go\l
\k{a}b-Schinzel type functional equation. \textsl{Abh. Math. Sem. Univ.
Hamburg,} \textbf{76} (2006), 91-98.

\noindent \lbrack Chu2] J. Chudziak, Semigroup-valued solutions of some
composite equations. \textsl{Aequat. Math.} \textbf{88} (2014), 183-198.%
\newline
\noindent \lbrack Cro] H. T. Croft, A question of limits, \textsl{Eureka}
\textbf{20} (1957), 11-13.\newline
\noindent \lbrack Dar1] G. Darboux, Sur la composition des forces en
statique. \textsl{Bull. Sci. Math.} \textbf{9} (1875), 281-299.\newline
\noindent \lbrack Dar2] G. Darboux, Sur le th\'{e}or\`{e}me fondamental de
la G\'{e}om\'{e}trie projective. \textsl{Math. Ann.} \textbf{17} (1880),
55-61.\newline
\noindent \lbrack FehLT] L. Feh\'{e}r, M. Laczkovich and G. Tardos, Croftian
sequences. \textsl{Acta Math. Hungar.} \textbf{56} (1990), 353--359.\newline
\noindent \lbrack FenN] J. E. Fenstad and D. Normann, On absolutely
measurable sets. \textsl{Fund. Math.} \textbf{81.2} (1973/74), 91--98.%
\newline
\noindent \lbrack GemH] D. Geman and J. Horowitz, Occupation densities.%
\textsl{\ Ann. Probab.} \textbf{8} (1980), 1-67.\newline
\noindent \lbrack Gro] {K.-G. Grosse-Erdmann, {An extension of the
Steinhaus-Weil theorem.} \textsl{Colloq. Math.} \textbf{57} (1989), 307--317.%
}\newline
\noindent \lbrack HulL] H. Hult and F. Linskog, Regular variation for
measures on metric spaces. \textsl{Publ. Inst. Math. Beograd} \textbf{90 (94)%
} (2006), 121-140.\newline
\noindent \lbrack HulLMS] H. Hult, F. Linskog, T. Mikosch and G.
Samorodnitsky, Functional large deviations for multivariate regularly
varying random walks. \textsl{Ann. Appl. Prob. }\textbf{15} (2005),
2651-2680.\newline
\noindent \lbrack Jab1] E. Jab\l o\'{n}ska, Continuous on rays solutions of
an equation of the Go\l \c{a}b-Schinzel type. \textsl{J. Math. Anal. Appl.}
\textbf{375} (2011), 223--229.\newline
\noindent \lbrack Jab2] E. Jab\l o\'{n}ska, On continuous solutions of an
equation of the Go\l \c{a}b-Schinzel type. \textsl{Bull. Aust. Math. Soc.}%
\textbf{\ 87} (2013), 10--17.\newline
\noindent \lbrack Jab3] E. Jab\l o\'{n}ska, On locally bounded above
solutions of an equation of the Go\l \c{a}b-Schinzel type. \textsl{Aequat.
Math.} \textbf{87} (2014),125--133.\newline
\noindent \lbrack Jab4] E. Jab\l o\'{n}ska, On continuous on rays solutions
of a composite-type equation. \textsl{Aequat. Math.} \textbf{89} (2015),
583--590.\newline
\noindent \lbrack Jab5] E. Jab\l o\'{n}ska, On solutions of some
generalizations of the Go\l \c{a}b-Schinzel equation. In: \textsl{Functional
equations in mathematical analysis}, ed. J. Brzd\k{e}k, Th. M. Rassias,
Springer Optim. Appl. \textbf{52} (2012), 509-521.\newline
\noindent \lbrack Jav] P. Javor, On the general solution of the functional
equation $f(x+yf(x))=f(x)f(y).$ \textsl{Aequat. Math.} \textbf{1} (1968),
235-238.\newline
\noindent \lbrack Kan] A. Kanamori, \textsl{The higher infinity. Large
cardinals in set theory from their beginnings}. Springer, 2$^{\text{nd}}$
ed. 2003 (1$^{\text{st}}$ ed. 1994).\newline
\noindent \lbrack Kec] A. S. Kechris, \textsl{Classical descriptive set
theory}. Grad. Texts in Math. \textbf{156}, Springer, 1996.\newline
\noindent \lbrack Ken] D. G. Kendall, Delphic semi-groups, infinitely
divisible regenerative phenomena, and the arithmetic of p-functions. \textsl{%
Z. Wahrscheinlich. verw. Geb.} \textbf{9} (1968), 163--195 (reprinted in
\textsl{Stochastic Analysis}: Rollo Davidson Memorial Volume (eds D. G.
Kendall, E. F. Harding), Wiley (1973), 73-114).\newline
\noindent \lbrack Kin1] J. F. C. Kingman, Ergodic properties of
continuous-time Markov processes and their discrete skeletons. \textsl{Proc.
London Math. Soc.} (3) 13 (1963), 593--604.\newline
\noindent \lbrack Kin2] J. F. C. Kingman, A note on limits of continuous
functions.\textsl{\ Quart. J. Math. Oxford Ser.} (2) \textbf{15} (1964),
279--282.\newline
\noindent \lbrack Kor] J. Korevaar, \textsl{Tauberian theorems: A century of
development}. Grundl. math. Wiss. \textbf{329}, Springer, 2004.\newline
\noindent \lbrack Kuc1] Marcin E. Kuczma, Differentiation of implicit
functions and Steinhaus' theorem in topological measure spaces. \textsl{%
Colloq. Math. }\textbf{39} (1978), 95--107, 189.\newline
\noindent \lbrack Kuc2] Marek Kuczma, Functional equations on restricted
domains. \textsl{Aequationes Math.} \textbf{18} (1978), 1--34.\newline
\noindent \lbrack Kuc3] Marek Kuczma, \textsl{An introduction to the theory
of functional equations and inequalities. Cauchy's equation and Jensen's
inequality.} 2nd ed., Birkh\"{a}user, 2009 [1st ed. PWN, Warszawa, 1985].%
\newline
\noindent \lbrack Mar] E. Marczewski, On translations of sets and a theorem
of Steinhaus. \textsl{Rocz. Pol. Tow. Mat.} (\textsl{Prace Mat.})\textbf{\
I.2} (1955), 256-263. (Transl. in: \textsl{Edward Marczewski -- Collected
Math. Papers}, 484-490, Pol. Acad. Sci, 1996.)\newline
\noindent \lbrack Mil] H. I. Miller. Generalization of a result of Borwein
and Ditor. \textsl{Proc. Amer. Math. Soc.} \textbf{105} (1989), 889--893.%
\newline
\noindent \lbrack NatO] T. Natkaniec and A. J. Ostaszewski, Embeddability of
countable sets by translation: interior point theorems, applications. In
preparation.\newline
\noindent \lbrack Ost1] A. J. Ostaszewski, Regular variation, topological
dynamics, and the Uniform Boundedness Theorem. \textsl{Top. Proc.}, \textbf{%
36} (2010), 305-336.\newline
\noindent \lbrack Ost2] A. J. Ostaszewski, Beurling regular variation, Bloom
dichotomy, and the Go\l \k{a}b-Schinzel functional equation. \textsl{Aequat.
Math.} \textbf{89} (2015), 725-744. \newline
\noindent \lbrack Ost3] A. J. Ostaszewski, Homomorphisms from Functional
Equations: The Goldie Equation. \textsl{Aequat. Math. }\textbf{90} (2016),
427-448 (arXiv: 1407.4089).\newline
\noindent \lbrack Ost4] A. J. Ostaszewski, Homomorphisms from Functional
Equations in Probability. In: \textsl{Developments in Functional Equations
and Related Topics}, ed. J. Brzd\k{e}k et al, Springer Optim. Appl. \textbf{%
124} (2017), 171-213.\newline
\noindent \lbrack Oxt] J. C. Oxtoby: \textsl{Measure and category.} 2nd ed.
Graduate Texts in Math. \textbf{2}, Springer, 1980.\newline
\noindent \lbrack Pet] {B. J. Pettis, {On continuity and openness of
homomorphisms in topological groups.} \textsl{Ann. of Math. }(2) \textbf{52}
(1950), 293--308.}\newline
\noindent \lbrack Pic] {S. Piccard, {Sur les ensembles de distances des
ensembles de points d'un espace Euclidien.\ }\textsl{M\'{e}m. Univ. Neuch%
\^{a}tel} \textbf{13}, 212 pp. 1939.}\newline
\noindent \lbrack Pop] C. G. Popa, Sur l'\'{e}quation fonctionelle $%
f[x+yf(x)]=f(x)f(y).$ \textsl{Ann. Polon. Math.} \textbf{17} (1965), 193-198.%
\newline
\noindent \lbrack Ree] D. Reem, Remarks on the Cauchy functional equation
and variations of it. \textsl{Aequationes Math.} \textbf{91} (2017),
237--264.\newline
\noindent \lbrack Rog] C. A. Rogers et al., \textsl{Analytic sets}. Academic
Press, 1980.\newline
\noindent \lbrack Sak] S. Saks,\textsl{\ Theory of the integra}l. Hafner,
New York, 1937 (translation of Monografie Matematyczne VII, 1937, 1$^{\text{%
st}}$ ed. Mono. Mat. II, 1933).

\noindent \lbrack Sen] E. Seneta, Sequential criteria for regular variation.
\textsl{Quart. J. Math. Oxford} Ser. (2) \textbf{22} (1971), 565--570.

\noindent \lbrack Sma] C. V. Smallwood, Approximate upper and lower limits.
\textsl{J. Math. Anal. App.} \textbf{37} (1972), 223-227.

\noindent \lbrack Sol] R. M. Solovay, A model of set-theory in which every
set of reals is Lebesgue measurable. \textsl{Ann. of Math.} (2) \textbf{92}
(1970), 1--56.\newline
\noindent \lbrack Vin] V. A. Vinokurov, On the limit at infinity of a
continuous function (in Russian). \textsl{Mat. Zametki }\textbf{1} (1967),
277-282.\newline
\noindent \lbrack Zin] R. E. Zink. On semicontinuous functions and Baire
functions. \textsl{Trans. Amer. Math. Soc.} \textbf{117} (1965), 1-6 .

\bigskip

\bigskip

\noindent Mathematics Department, Imperial College, London SW7 2AZ;
n.bingham@ic.ac.uk \newline
Mathematics Department, London School of Economics, Houghton Street, London
WC2A 2AE; A.J.Ostaszewski@lse.ac.uk\newpage

\end{document}